\documentclass[a4paper,12pt]{article} 
\usepackage{amssymb}
\usepackage{amscd}
\usepackage{amsmath} 
\usepackage{graphicx}
\usepackage{latexsym} 
\usepackage{enumerate}

\usepackage{theorem}
\theoremstyle{plain}
\theorembodyfont{\itshape}
\newtheorem{theorem}{Theorem}[section]
\newtheorem{proposition}[theorem]{Proposition}
\newtheorem{lemma}[theorem]{Lemma}

\theorembodyfont{\rmfamily}

\newtheorem{remark}[theorem]{Remark}

\def\Aff{\mathrm{Aff}}
\def\Aut{\mathrm{Aut}}

\def\det{\mathrm{det}}
\def\End{\mathrm{End}}

\def\PVI{\mathrm{P}_{\mathrm{VI}}}
\def\Per{\mathrm{Per}} 
\def\PS{\mathrm{PS}}

\def\Res{\mathrm{Res}}
\def\RH{\mathrm{RH}} 
\def\rh{\mathrm{rh}}

\def\Tr{\mathrm{Tr}}

\def\E{\mathcal{E}}

\def\K{\mathcal{K}}

\def\M{\mathcal{M}}
\def\O{\mathcal{O}}
\def\R{\mathcal{R}}

\def\Sol{\mathcal{S}}

\def\C{\mathbb{C}}

\def\N{\mathbb{N}}
\def\P{\mathbb{P}}

\def\bR{\mathbb{R}}
\def\Z{\mathbb{Z}}
\def\Wall{\mathbf{Wall}} 

\def\b{\beta}
\def\ga{\gamma}
\def\k{\kappa}
\def\l{\lambda}
\def\si{\sigma}
\def\th{\theta}
\def\vD{\varDelta}
\def\vG{\varGamma}

\def\Th{\Theta}
\def\ve{\varepsilon}
\def\carl{\circlearrowleft}
\def\car{\curvearrowright}

\def\ra{\rightarrow}

\def\ol{\overline}

\def\ci{\circ}
\def\dfrac#1#2{{\displaystyle\frac{#1}{#2}}}

\def\la{\langle}
\def\ra{\rangle}

\def\-{\phantom{-}}
\title{\bf Periodic Solutions to Painlev\'e VI and \\ 
Dynamical System on Cubic Surface} 
\author{Katsunori Iwasaki and Takato Uehara \\ \\
Graduate School of Mathematics, Kyushu University \\
6-10-1 Hakozaki, Higashi-ku, Fukuoka 812-8581 
Japan\thanks{E-mail addresses: {\tt iwasaki@math.kyushu-u.ac.jp} \ 
and \ {\tt ma205003@math.kyushu-u.ac.jp}}} 
\date{Dedicated to Professor Masuo Hukuhara on his $100$th birthday} 
\begin{document}
\maketitle
\begin{abstract} 
The number of periodic solutions to Painlev\'e VI along a 
Pochhammer loop is counted exactly. 
It is shown that the number grows exponentially with period, 
where the growth rate is determined explicitly. 
Principal ingredients of the computation are a moduli-theoretical 
formulation of Painlev\'e VI, a Riemann-Hilbert correspondence, 
the dynamical system of a birational map on a cubic surface, and 
the Lefschetz fixed point formula. 
\end{abstract} 
\section{Introduction}  \label{sec:intro} 
Painlev\'e equations and dynamical systems on complex surfaces 
are two subjects of mathematics which have been investigated 
actively in recent years. 
In this paper we shall demonstrate a substantial relation between 
them by presenting a fruitful application to the former subject 
of the latter. 
We begin with stating our motivation on the side 
of Painlev\'e equations. 
\par 
The global structure of the sixth Painlev\'e equation $\PVI(\k)$, 
especially the multivalued character of its solutions 
is an important issue in the study of Painlev\'e 
equations. 
In this respect, several authors 
\cite{Boalch1,Boalch2,DM,Hitchin1,Hitchin2,Kitaev,Mazzocco} have 
been interested in finding algebraic solutions, because they 
offer a simplest class of solutions with clear global structure 
in the sense that they have only finitely many branches under 
analytic continuations along all loops in the domain 
\begin{equation} \label{eqn:X}
X = \P^1-\{0,1,\infty\}. 
\end{equation}
\par 
In another direction of promising research, we are 
interested in periodic solutions along a single loop, namely, 
in those solutions which are finitely many-valued along a single 
loop chosen particularly. 
Given such a loop, we shall discuss the following problems: 
\begin{itemize} 
\item How many solutions can be periodic of period $N$ 
among all solutions to $\PVI(\k)$? 
\vspace{-0.2cm} 
\item How rapidly does that number grow as the period 
$N$ tends to infinity?
\end{itemize}
\par
For such a loop we take a Pochhammer loop $\wp$ as in 
Figure \ref{fig:pochhammer}. 
If $\ell_0$ and $\ell_1$ are standard generators of $\pi_1(X,x)$ 
as in Figure \ref{fig:loop}, then $\wp$ is a loop homotopic to 
the commutator 
\[
[\ell_0, \ell_1^{-1}] = \ell_0 \ell_1^{-1} \ell_0^{-1} \ell_1. 
\]
It is a typical loop which often appears in mathematics due to 
the property that any abelian representation of $\pi_1(X,x)$ 
is killed along this loop. 
For example, it is used as an integration contour of Euler 
integral representation of hypergeometric 
functions \cite{IKSY}. 
In the context of this article the Pochhammer loop will be 
closely related to a certain birational map of a cubic surface 
whose dynamics is quite relevant to understanding the 
global structure of the sixth Painlev\'e equation 
(see discussions in Sections \ref{sec:involution} 
and \ref{sec:dynamics}). 
\par 
We remark that the same problem for a simplest loop, namely 
for a loop $\ell_0$ or $\ell_1$ in Figure \ref{fig:loop} or a 
loop $\ell_{\infty} = (\ell_0\ell_1)^{-1}$ around the point at 
infinity, is not interesting or even meaningless, 
because for any $N > 1$, there are infinitely many periodic 
solutions of period $N$ along it. 
In fact it turns out that they are parametrized by points on 
certain complex curves and hence their 
cardinality is that of a continuum \cite{Iwasaki5}. 
On the contrary, along the Pochhammer loop $\wp$, 
the cardinality of the periodic solutions of period $N$ turns 
out to be finite for every $N \in \N := \{1,2,3,\dots\}$ and 
hence our problem certainly makes sense. 
\begin{figure}[t]
\begin{center}
\unitlength 0.1in
\begin{picture}(42.15,18.36)(4.75,-19.21)
%
\special{pn 20}%
\special{sh 0.600}%
\special{ar 1390 1000 67 67  0.0000000 6.2831853}%
%
\special{pn 20}%
\special{sh 0.600}%
\special{ar 3780 1000 67 67  0.0000000 6.2831853}%
\put(13.1000,-13.0000){\makebox(0,0)[lb]{$0$}}%
\put(37.2000,-13.0000){\makebox(0,0)[lb]{$1$}}%
%
\special{pn 20}%
\special{pa 1980 1100}%
\special{pa 2910 1100}%
\special{fp}%
%
\special{pn 20}%
\special{pa 1420 400}%
\special{pa 1340 400}%
\special{fp}%
\special{sh 1}%
\special{pa 1340 400}%
\special{pa 1407 420}%
\special{pa 1393 400}%
\special{pa 1407 380}%
\special{pa 1340 400}%
\special{fp}%
%
\special{pn 20}%
\special{pa 790 960}%
\special{pa 790 1050}%
\special{fp}%
\special{sh 1}%
\special{pa 790 1050}%
\special{pa 810 983}%
\special{pa 790 997}%
\special{pa 770 983}%
\special{pa 790 1050}%
\special{fp}%
%
\special{pn 20}%
\special{pa 1350 1610}%
\special{pa 1440 1610}%
\special{fp}%
\special{sh 1}%
\special{pa 1440 1610}%
\special{pa 1373 1590}%
\special{pa 1387 1610}%
\special{pa 1373 1630}%
\special{pa 1440 1610}%
\special{fp}%
%
\special{pn 20}%
\special{pa 2475 1100}%
\special{pa 2536 1100}%
\special{fp}%
\special{sh 1}%
\special{pa 2536 1100}%
\special{pa 2469 1080}%
\special{pa 2483 1100}%
\special{pa 2469 1120}%
\special{pa 2536 1100}%
\special{fp}%
%
\special{pn 20}%
\special{pa 3740 1620}%
\special{pa 3850 1610}%
\special{fp}%
\special{sh 1}%
\special{pa 3850 1610}%
\special{pa 3782 1596}%
\special{pa 3797 1615}%
\special{pa 3785 1636}%
\special{pa 3850 1610}%
\special{fp}%
%
\special{pn 20}%
\special{pa 4410 1040}%
\special{pa 4410 960}%
\special{fp}%
\special{sh 1}%
\special{pa 4410 960}%
\special{pa 4390 1027}%
\special{pa 4410 1013}%
\special{pa 4430 1027}%
\special{pa 4410 960}%
\special{fp}%
%
\special{pn 20}%
\special{pa 3830 390}%
\special{pa 3770 390}%
\special{fp}%
\special{sh 1}%
\special{pa 3770 390}%
\special{pa 3837 410}%
\special{pa 3823 390}%
\special{pa 3837 370}%
\special{pa 3770 390}%
\special{fp}%
%
\special{pn 20}%
\special{pa 1430 1920}%
\special{pa 1360 1910}%
\special{fp}%
\special{sh 1}%
\special{pa 1360 1910}%
\special{pa 1423 1939}%
\special{pa 1413 1918}%
\special{pa 1429 1900}%
\special{pa 1360 1910}%
\special{fp}%
%
\special{pn 20}%
\special{pa 480 1060}%
\special{pa 480 960}%
\special{fp}%
\special{sh 1}%
\special{pa 480 960}%
\special{pa 460 1027}%
\special{pa 480 1013}%
\special{pa 500 1027}%
\special{pa 480 960}%
\special{fp}%
%
\special{pn 20}%
\special{pa 1360 90}%
\special{pa 1440 90}%
\special{fp}%
\special{sh 1}%
\special{pa 1440 90}%
\special{pa 1373 70}%
\special{pa 1387 90}%
\special{pa 1373 110}%
\special{pa 1440 90}%
\special{fp}%
%
\special{pn 20}%
\special{pa 3750 120}%
\special{pa 3860 110}%
\special{fp}%
\special{sh 1}%
\special{pa 3860 110}%
\special{pa 3792 96}%
\special{pa 3807 115}%
\special{pa 3795 136}%
\special{pa 3860 110}%
\special{fp}%
%
\special{pn 20}%
\special{pa 4680 950}%
\special{pa 4690 1060}%
\special{fp}%
\special{sh 1}%
\special{pa 4690 1060}%
\special{pa 4704 992}%
\special{pa 4685 1007}%
\special{pa 4664 995}%
\special{pa 4690 1060}%
\special{fp}%
%
\special{pn 20}%
\special{pa 3830 1900}%
\special{pa 3740 1900}%
\special{fp}%
\special{sh 1}%
\special{pa 3740 1900}%
\special{pa 3807 1920}%
\special{pa 3793 1900}%
\special{pa 3807 1880}%
\special{pa 3740 1900}%
\special{fp}%
\put(24.9000,-15.9000){\makebox(0,0)[lb]{$\wp$}}%
%
\special{pn 20}%
\special{ar 3790 1010 894 894  3.0396743 6.2831853}%
\special{ar 3790 1010 894 894  0.0000000 2.9289477}%
%
\special{pn 20}%
\special{pa 2280 1190}%
\special{pa 2910 1190}%
\special{fp}%
%
\special{pn 20}%
\special{pa 1980 910}%
\special{pa 3200 910}%
\special{fp}%
%
\special{pn 20}%
\special{ar 1390 1010 598 598  0.1513754 6.1152894}%
%
\special{pn 20}%
\special{ar 1390 1010 911 911  0.1995554 6.2831853}%
%
\special{pn 20}%
\special{pa 2310 1010}%
\special{pa 3190 1010}%
\special{fp}%
%
\special{pn 20}%
\special{ar 3800 1010 620 620  3.3067413 6.2831853}%
\special{ar 3800 1010 620 620  0.0000000 3.1415927}%
%
\special{pn 20}%
\special{pa 2710 910}%
\special{pa 2650 910}%
\special{fp}%
\special{sh 1}%
\special{pa 2650 910}%
\special{pa 2717 930}%
\special{pa 2703 910}%
\special{pa 2717 890}%
\special{pa 2650 910}%
\special{fp}%
%
\special{pn 20}%
\special{pa 2720 1190}%
\special{pa 2650 1190}%
\special{fp}%
\special{sh 1}%
\special{pa 2650 1190}%
\special{pa 2717 1210}%
\special{pa 2703 1190}%
\special{pa 2717 1170}%
\special{pa 2650 1190}%
\special{fp}%
%
\special{pn 20}%
\special{pa 2470 1010}%
\special{pa 2540 1000}%
\special{fp}%
\special{sh 1}%
\special{pa 2540 1000}%
\special{pa 2471 990}%
\special{pa 2487 1008}%
\special{pa 2477 1029}%
\special{pa 2540 1000}%
\special{fp}%
\end{picture}%
\end{center}
\caption{Pochhammer loop $\wp$} 
\label{fig:pochhammer} 
\end{figure}
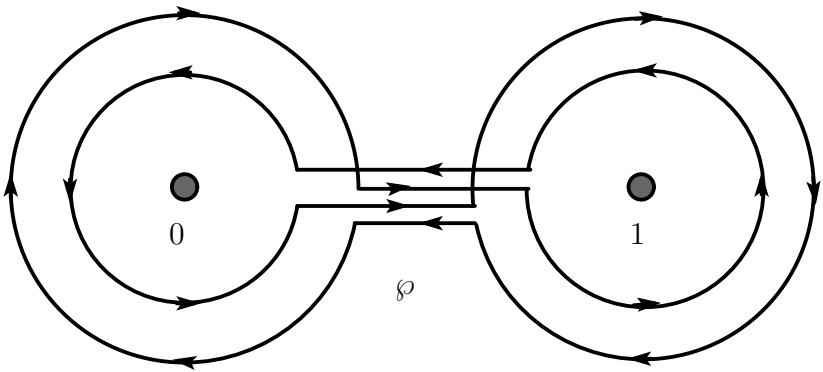
\begin{figure}[b]
\begin{center}
\unitlength 0.1in
\begin{picture}(36.05,13.50)(7.80,-17.50)
%
\special{pn 20}%
\special{sh 0.600}%
\special{ar 1390 1000 67 67  0.0000000 6.2831853}%
%
\special{pn 20}%
\special{sh 0.600}%
\special{ar 2600 1000 67 67  0.0000000 6.2831853}%
%
\special{pn 20}%
\special{sh 0.600}%
\special{ar 3780 1000 67 67  0.0000000 6.2831853}%
%
\special{pn 20}%
\special{ar 1390 1010 610 610  0.0782711 6.1851400}%
%
\special{pn 20}%
\special{ar 3780 1010 600 600  3.2429392 6.2831853}%
\special{ar 3780 1010 600 600  0.0000000 3.0750245}%
%
\special{pn 20}%
\special{pa 2000 950}%
\special{pa 3190 950}%
\special{fp}%
%
\special{pn 20}%
\special{pa 2000 1050}%
\special{pa 3190 1050}%
\special{fp}%
%
\special{pn 20}%
\special{pa 780 980}%
\special{pa 790 1070}%
\special{fp}%
\special{sh 1}%
\special{pa 790 1070}%
\special{pa 803 1002}%
\special{pa 784 1017}%
\special{pa 763 1006}%
\special{pa 790 1070}%
\special{fp}%
%
\special{pn 20}%
\special{pa 4380 1040}%
\special{pa 4380 960}%
\special{fp}%
\special{sh 1}%
\special{pa 4380 960}%
\special{pa 4360 1027}%
\special{pa 4380 1013}%
\special{pa 4400 1027}%
\special{pa 4380 960}%
\special{fp}%
\put(13.1000,-13.0000){\makebox(0,0)[lb]{$0$}}%
\put(37.2000,-13.0000){\makebox(0,0)[lb]{$1$}}%
\put(25.3000,-13.0000){\makebox(0,0)[lb]{$x$}}%
\put(13.0000,-19.2000){\makebox(0,0)[lb]{$\ell_0$}}%
\put(37.0000,-19.1000){\makebox(0,0)[lb]{$\ell_1$}}%
%
\special{pn 20}%
\special{pa 1440 410}%
\special{pa 1350 400}%
\special{fp}%
\special{sh 1}%
\special{pa 1350 400}%
\special{pa 1414 427}%
\special{pa 1403 406}%
\special{pa 1418 387}%
\special{pa 1350 400}%
\special{fp}%
%
\special{pn 20}%
\special{pa 1310 1620}%
\special{pa 1430 1620}%
\special{fp}%
\special{sh 1}%
\special{pa 1430 1620}%
\special{pa 1363 1600}%
\special{pa 1377 1620}%
\special{pa 1363 1640}%
\special{pa 1430 1620}%
\special{fp}%
%
\special{pn 20}%
\special{pa 3740 1620}%
\special{pa 3800 1610}%
\special{fp}%
\special{sh 1}%
\special{pa 3800 1610}%
\special{pa 3731 1601}%
\special{pa 3747 1619}%
\special{pa 3738 1641}%
\special{pa 3800 1610}%
\special{fp}%
%
\special{pn 20}%
\special{pa 3830 410}%
\special{pa 3760 410}%
\special{fp}%
\special{sh 1}%
\special{pa 3760 410}%
\special{pa 3827 430}%
\special{pa 3813 410}%
\special{pa 3827 390}%
\special{pa 3760 410}%
\special{fp}%
\end{picture}%
\end{center}
\caption{Standard generators of $\pi_1(X,x)$} 
\label{fig:loop}
\end{figure}
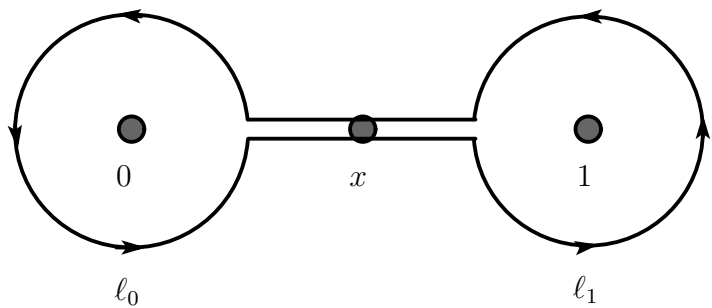
\par 
In this article we shall exactly count the number of periodic 
solutions to $\PVI(\k)$ along the Pochhammer loop $\wp$ under 
a certain generic assumption on the parameters $\k$. 
In particular we shall show that the number grows exponentially 
as the period tends to infinity, with the growth rate 
determined explicitly (see Theorem \ref{thm:main}). 
As is already mentioned, this result is the fruits of a good 
application to Painlev\'e equations of a dynamical system 
theory on complex surfaces as developed in \cite{DF,Favre}. 
Algebraic geometry of Painlev\'e equations, especially 
a moduli-theoretical interpretation of Painlev\'e dynamics 
\cite{IIS1,IIS2} is also an essential ingredient of our 
work. 
\par 
After stating the main result of this article in 
Section \ref{sec:result}, we shall develop the story 
of this article in the following manner.  
First, $\PVI(\k)$ is formulated as a flow, {\sl Painlev\'e flow}, 
on a moduli space of stable parabolic connections 
(Section \ref{sec:mspc}). 
Secondly, it is conjugated to an isomonodromic flow 
on a moduli space of monodromy representations via a 
Riemann-Hilbert correspondence (Section \ref{sec:RHC}). 
Thirdly, with a natural identification of the representation 
space with a cubic surface, the Poincar\'e section of $\PVI(\k)$ 
is conjugated to the dynamical system of a group action on the 
cubics (Section \ref{sec:cubic}). 
Especially, analytic continuation along the Pochhammer 
loop is connected with a distinguished transformation, 
called a {\sl `Coxeter' transformation}, of the group action. 
Fourthly, main properties of our dynamical system on the 
cubics are established from the standpoint of birational 
surface dynamics. 
Fifthly, the number of periodic points of the Coxeter 
transformation is counted by using the Lefschetz 
fixed point formula (Section \ref{sec:lefschetz}). 
Then, back to the original phase space of $\PVI(\k)$, 
we arrive at our final goal, that is, the exact number of 
periodic solutions to $\PVI(\k)$ of any period along the 
Pochhammer loop $\wp$. 
\par 
The authors would be happy if this article could give a new 
insight into the global structure of the sixth Painlev\'e equation. 
They are grateful to Yutaka Ishii for valuable discussions. 
\section{Main Result} \label{sec:result}
Let us describe our main result in more detail. 
To this end we recall that the sixth Painlev\'e equation 
$\PVI(\k)$ in its Hamiltonian form is a system of nonlinear 
differential equations 
\begin{equation} \label{eqn:PVI}
\dfrac{d q}{d x} = \dfrac{\partial H(\k)}{\partial p}, 
\qquad 
\dfrac{d p}{d x} = -\dfrac{\partial H(\k)}{\partial q}, 
\end{equation}
with an independent variable $x \in X$ and unknown functions 
$(q(x), p(x))$, depending on complex parameters 
$\k = (\k_0,\k_1,\k_2,\k_3,\k_4)$ in a four-dimensional affine 
space 
\[
\K := \{\, \k = (\k_0,\k_1,\k_2,\k_3,\k_4) \in \C^5 \,:\, 
2 \k_0 + \k_1 + \k_2 + \k_3 +\k_4 = 1 \,\},  
\]
where the Hamiltonian $H(\k) = H(q,p,x;\k)$ is given by 
\[
\begin{array}{rcl}
x(x-1) H(\k) &=&  
(q_0q_1q_x) p^2 - \{\k_1q_1q_x + (\k_2-1)q_0q_1 + \k_3q_0q_x\} p 
+ \k_0(\k_0+\k_4) q_x, 
\end{array}
\]
with $q_{\nu} = q - \nu$ for $\nu \in \{0, 1, x \}$. 
It is known that system (\ref{eqn:PVI}) enjoys the Painlev\'e 
property, that is, any meromorphic solution germ at a base point 
$x \in X$ of system (\ref{eqn:PVI}) admits a global analytic 
continuation along any path emanating from $x$ as a meromorphic 
function. 
\par 
Geometrically, the sixth Painlev\'e equation $\PVI(\k)$ is 
formulated as a holomorphic uniform foliation on the total 
space of a fibration of certain smooth, quasi-projective, 
rational surfaces, 
\begin{equation} \label{eqn:phase}
\pi_{\k} : M(\k) \to X, 
\end{equation} 
transversal to each fiber of the fibration. 
We refer to \cite{AL,IIS1,IIS2,Okamoto,STT,Sakai} for the 
detailed accounts of the space $M(\k)$. 
Especially the papers \cite{IIS1,IIS2} give a comprehensive 
description of it as a moduli space of stable parabolic 
connections. 
The fiber $M_x(\k)$ over $x \in X$, called the 
space of initial conditions at time $x$, parametrizes all 
the solution germs at $x$ of equation (\ref{eqn:PVI}) most 
precisely, completing the na\"{\i}ve and incomplete space 
$\C^2$ of initial values $(q,p)$ at the points $x$. 
Given a loop $\ga \in \pi_1(X,x)$, the horizontal lifts 
of the loop $\ga$ along the foliation induces  
a biholomorphism $\ga_* : M_x(\k) \to M_x(\k)$, 
called the Poincar\'e return map along $\ga$, which 
depends only on the homotopy class of $\ga$. 
Then the global structure of the sixth Painlev\'e equation 
$\PVI(\k)$ is described by the group homomorphism 
\begin{equation} \label{eqn:PSx}
\PS_x(\k) \,\, : \,\, \pi_1(X,x) \to \Aut\, M_x(\k), 
\qquad \ga \mapsto \ga_*, 
\end{equation}
which is referred to as the {\sl Poincar\'e section} of the 
sixth Painlev\'e dynamics $\PVI(\k)$. 
\par 
In this article we are interested in analytic continuations 
of solutions to equation (\ref{eqn:PVI}) along the Pochhammer 
loop $\wp$, namely, in the iteration of the Poincar\'e return 
map $\wp_*$ along $\wp$. 
The Poincar\'e return map along the Pochhammer loop is referred 
to as the {\sl Pochhammer-Poincar\'e map}. 
Given any $N \in \N$, let $\Per_N(\k)$ be the set of 
all initial points $Q \in M_x(\k)$ that come back 
to the original positions after the $N$-th iterate of the 
Pochhammer-Poincar\'e map $\wp_*$,  
\begin{equation} \label{eqn:per}
\Per_N(\k) := \{\, Q \in M_x(\k) \,:\, 
\wp_*^N(Q) = Q \, \}. 
\end{equation}
The aim of this article is to count the number of $\Per_N(\k)$ 
and to find out its growth rate as the period 
$N$ tends to infinity. 
\par
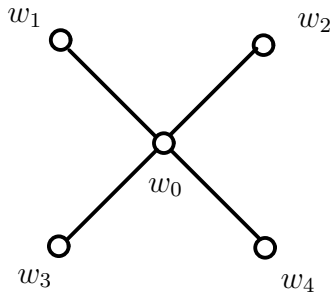
\begin{figure}[t]
\begin{center}
\unitlength 0.1in
\begin{picture}(14.83,13.55)(8.95,-14.45)
%
\special{pn 20}%
\special{ar 1795 797 52 52  6.1412883 6.2831853}%
\special{ar 1795 797 52 52  0.0000000 6.1180366}%
%
\special{pn 20}%
\special{ar 1255 257 52 52  6.1412883 6.2831853}%
\special{ar 1255 257 52 52  0.0000000 6.1180366}%
%
\special{pn 20}%
\special{ar 2317 284 52 52  6.1412883 6.2831853}%
\special{ar 2317 284 52 52  0.0000000 6.1180366}%
%
\special{pn 20}%
\special{ar 1246 1337 52 52  6.1412883 6.2831853}%
\special{ar 1246 1337 52 52  0.0000000 6.1180366}%
%
\special{pn 20}%
\special{ar 2326 1346 52 52  6.1412883 6.2831853}%
\special{ar 2326 1346 52 52  0.0000000 6.1180366}%
%
\special{pn 20}%
\special{pa 1300 311}%
\special{pa 1750 761}%
\special{fp}%
%
\special{pn 20}%
\special{pa 1840 842}%
\special{pa 2299 1301}%
\special{fp}%
%
\special{pn 20}%
\special{pa 2281 302}%
\special{pa 1831 752}%
\special{fp}%
%
\special{pn 20}%
\special{pa 1759 833}%
\special{pa 1291 1301}%
\special{fp}%
\put(18.0400,-10.2600){\makebox(0,0){$w_0$}}%
\put(26.7700,-1.1700){\makebox(0,0)[rt]{$w_2$}}%
\put(24.9700,-15.3000){\makebox(0,0){$w_4$}}%
\put(11.2000,-15.0300){\makebox(0,0){$w_3$}}%
\put(9.7600,-0.9000){\makebox(0,0)[lt]{$w_1$}}%
\end{picture}%
\end{center}
\caption{Dynkin diagram of type $D_4^{(1)}$}
\label{fig:dynkin} 
\end{figure}

 
To avoid certain technical difficulties 
(see Remark \ref{rem:main2}), we make a generic 
assumption on the parameters $\k \in \K$. 
To this end we recall an affine Weyl group structure 
of the parameter space \cite{IIS1,Iwasaki4}. 
The affine space $\K$ is identified with the linear 
space $\C^4$ by the isomorphism 
\[
\K \to \C^4, \quad \k = (\k_0,\k_1,\k_2,\k_3,\k_4) 
\mapsto (\k_1,\k_2,\k_3,\k_4), 
\]
where the latter space $\C^4$ is equiped with the 
standard (complex) Euclidean inner product. 
For each $i \in \{0,1,2,3,4\}$, let $w_i : \K \to \K$ be 
the orthogonal reflection having $\{\k_i =0\}$ as its 
reflecting hyperplane, with respect to the inner 
product mentioned above. 
Then the group generated by $w_0$, $w_1$, $w_2$, $w_3$, 
$w_4$ is an affine Weyl group of type $D_4^{(1)}$, 
\[
W(D_4^{(1)}) = \la w_0, w_1, w_2, w_3, w_4 \ra 
\car \K.   
\]
corresponding to the Dynkin diagram in Figure 
\ref{fig:dynkin}. 
The reflecting hyperplanes of all reflections in the 
group $W(D_4^{(1)})$ are given by affine linear relations 
\[
\k_i = m, \qquad \k_1 \pm \k_2 \pm \k_3 \pm \k_4 = 2m+1 
\qquad (i \in \{1,2,3,4\}, \, m \in \Z), 
\]
with any choice of signs $\pm$. 
Let $\Wall$ be the union of all those hyperplanes. 
Then the generic condition to be imposed on parameters is 
that $\k$ should lie outside $\Wall$; this is a necessary 
and sufficient condition for $\PVI(\k)$ to admit no 
Riccati solutions \cite{IIS1}. 
\par 
Now the main theorem of this article is stated as follows. 
\begin{theorem} \label{thm:main} 
For any $\k \in \K - \Wall$ the cardinality of the set 
$\Per_N(\k)$ is given by 
\begin{equation} \label{eqn:main}
\mathrm{\#} \, \Per_N(\k) = 
(9+4\sqrt{5})^N +(9-4\sqrt{5})^N + 4 
\qquad (N \in \N). 
\end{equation}
\end{theorem} 
\begin{remark} \label{rem:main1} 
It should be noted that formula (\ref{eqn:main}) is 
rewritten as 
\[
\mathrm{\#}    \, \Per_N(\k) - \{(9+4\sqrt{5})^N + 4\} =
(9+4\sqrt{5})^{-N}, 
\]
which means that the geometric sequence $(9+4\sqrt{5})^N$ 
shifted by $4$ approximates the cardinality of 
$\Per_N(\k)$ up to an exponentially decaying error 
term $(9+4\sqrt{5})^{-N}$, where the growth rate of 
cardinality and the decay rate of error term are 
given by the same number $9+4\sqrt{5}$. 
Moreover, since $9 \pm 4\sqrt{5}$ are the root of the 
quadratic equation $\l^2 - 18 \l + 1 = 0$, the formula 
(\ref{eqn:main}) is expressed as 
$\mathrm{\#} \, \Per_N(\k) = C_N + 4$, where the 
sequence $\{C_N\}$ is defined recursively by 
\[
C_0 = 2, \qquad C_1 = 18, \qquad 
C_{N+2} - 18 \, C_{N+1} + C_N = 0. 
\]
\end{remark}
\begin{remark} \label{rem:main3} 
Our main theorem can also be stated in terms of a dynamical zeta 
function. 
Indeed, as a generating expression of formula (\ref{eqn:main}) 
for all $N \in \N$, we have 
\[ 
Z_{\k}(z) := \exp\left(\sum_{N=1}^{\infty} 
\mathrm{\#}    \, \Per_N(\k) \, \dfrac{z^N}{N} \right) 
= \dfrac{1}{(1-z)^4(1 - 18 z + z^2)}. 
\]
\end{remark}
\begin{remark} \label{rem:main2} 
In this article we restrict our attention to the 
generic case $\k \in \K-\Wall$ only, leaving 
the nongeneric case $\k \in \Wall$ untouched. 
The difference between the generic case and the nongeneric 
case lies in the fact that the Riemann-Hilbert correspondence 
to be used in the proof becomes a biholomorphism in the former 
case, while it gives an analytic minimal resolution of Klein 
singularities in the latter case (see Remark \ref{rem:RH}). 
The presence of singularities would make the treatment of 
the nongeneric case more complicated. 
However it is expected that the basic strategy developed in this 
article will also be effective in the nongeneric case. 
The relevant discussion will be made in another place. 
\end{remark}
\section{Moduli Space of Stable Parabolic Connections} 
\label{sec:mspc}
In order to describe the fibration (\ref{eqn:phase}), we first 
construct an auxiliary fibration 
$\pi_{\k} : \M(\k) \to T$ over the configuration 
space of mutually distinct, ordered, three points in $\C$, 
\[
T = \{\, t=(t_1,t_2,t_3) \in \C^3\,:\, t_i \neq t_j \,\, 
\text{for} \,\, i \neq j \, \},  
\]
and then reduce it to the original fibration (\ref{eqn:phase}). 
We put the fourth point $t_4$ at infinity. 
Given any $(t,\k) \in T \times \K$, a $(t,\k)$-parabolic 
connection is a quadruple $Q = (E,\nabla,\psi,l)$ such that
\begin{enumerate} 
\item $E$ is a rank $2$ vector bundle of degree $-1$ over $\P^1$, 
\item $\nabla : E \to E \otimes \Omega^1_{\P^1}(D_t)$ is a Fuchsian 
connection with pole divisor $D_t = t_1 + t_2 + t_3 + t_4$ and 
Riemann scheme as in Table \ref{tab:riemann}, where 
$t_4 = \infty$ as mentioned above, 
\item $\psi : \det \, E \to \O_{\P^1}(-t_4)$ is a horizontal 
isomorphism called a determinantal structure, where 
$\O_{\P^1}(-t_4)$ is equiped with the connection induced from 
$d : \O_{\P^1} \to \Omega^1_{\P^1}$,  
\item $l = (l_1,l_2,l_3,l_4)$ is a parabolic structure, namely, 
$l_i$ is an eigenline of $\Res_{t_i}(\nabla) \in \End(E_{t_i})$ 
corresponding to eigenvalue $\l_i$ (whose minus is the first 
exponent $-\l_i$ in Table \ref{tab:riemann}). 
\end{enumerate} 
\begin{table}[t]
\begin{center} 
\begin{tabular}{|c||c|c|c|c|}
\hline
\vspace{-0.3cm} & & & & \\
singularities & $t_1$ & $t_2$ & $t_3$ & $t_4$ \\[2mm]
\hline
\vspace{-0.3cm} & & & & \\
first exponent & $-\l_1$ & $-\l_2$ & $-\l_3$ & $-\l_4$ \\[2mm]
\hline
\vspace{-0.3cm} & & & & \\
second exponent & $\l_1$ & $\l_2$  & $\l_3$ & $\l_4-1$ \\[2mm]
\hline
\vspace{-0.3cm} & & & & \\
difference  & $\k_1$ & $\k_2$ & $\k_3$ & $\k_4$ \\[2mm]
\hline
\end{tabular}
\end{center}
\caption{Riemann scheme: $\k_i$ is the difference of the second 
exponent from the first.} 
\label{tab:riemann}
\end{table}
There exists a concept of stability for parabolic connections, 
with which the geometric invariant theory \cite{Mumford} can 
be worked out to establish the following theorem \cite{IIS1,IIS2}.  
\begin{theorem} \label{thm:moduli} 
For any $(t,\k) \in T \times \K$ there exists a fine moduli scheme 
$\M_t(\k)$ of stable $(t,\k)$-parabolic connections.  
The moduli space $\M_t(\k)$ is a smooth, irreducible, 
quasi-projective surface. 
As a relative setting over $T$, for any $\k \in \K$, there exists 
a family of moduli spaces 
\begin{equation} \label{eqn:family} 
\pi_{\k} : \M(\k) \rightarrow T 
\end{equation}
such that the projection $\pi_{\k}$ is a smooth morphism 
with fiber $\M_t(\k)$ over $t \in T$. 
\end{theorem} 
\par
Now the fibration (\ref{eqn:phase}) is defined to be the pull-back 
of (\ref{eqn:family}) by an injection 
\[
\iota : X \hookrightarrow T, \quad x \mapsto (0,x,1),  
\]
The group $\Aff(\C)$ of affine linear transformations on 
$\C$ acts diagonally on the configuration space $T$ 
and the quotient space $T/\Aff(\C)$ is isomorphic 
to $X$, with the quotient map given by 
\begin{equation} \label{eqn:reduction} 
r : T \to X, \quad t = (t_1,t_2,t_3) \mapsto 
x = \dfrac{t_2-t_1}{t_3-t_1}. 
\end{equation}
The map $r$ yields a trivial $\Aff(\C)$-bundle structure of $T$ 
over $X$ and the fibration (\ref{eqn:family}) is in turn the 
pull-back of the fibration (\ref{eqn:phase}) by the map $r$. 
Hence we have a commutative diagram 
\begin{equation} \label{cd:reduction} 
\begin{CD}
\M(\k) @>>> M(\k) \\
 @V \pi_{\k} VV  @VV \pi_{\k} V \\
T @>> r > X. 
\end{CD}
\end{equation}
In \cite{IIS1,IIS2} the Painlev\'e dynamics $\PVI(\k)$ is 
formulated as a holomorphic uniform foliations on the fibration 
(\ref{eqn:family}) which is compatible with the diagram 
(\ref{cd:reduction}). 
Thus the Poincar\'e section (\ref{eqn:PSx}) is reformulated as a 
group homomorphism
\begin{equation} \label{eqn:PSt} 
\PS_t(\k) \,\,:\,\, \pi_1(T,t) \to \Aut\,\M_t(\k). 
\end{equation}
\par
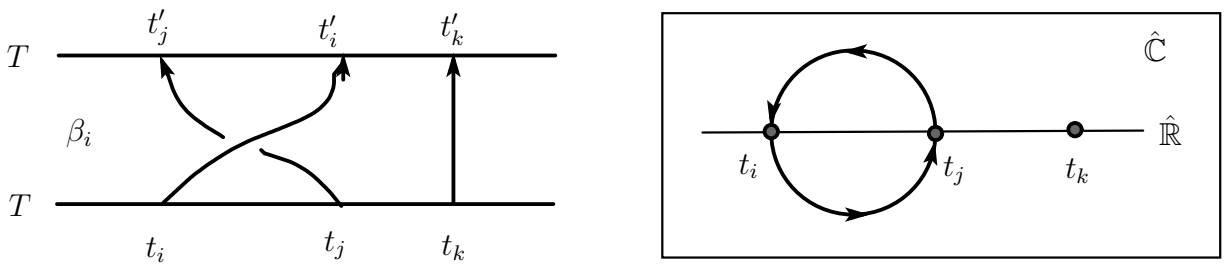
\begin{figure}[t]
\begin{center}
\unitlength 0.1in
\begin{picture}(63.50,12.85)(1.20,-14.60)
%
\special{pn 20}%
\special{pa 391 403}%
\special{pa 2990 403}%
\special{fp}%
%
\special{pn 20}%
\special{pa 383 1180}%
\special{pa 2981 1180}%
\special{fp}%
%
\special{pn 20}%
\special{pa 2457 1170}%
\special{pa 2457 423}%
\special{fp}%
\special{sh 1}%
\special{pa 2457 423}%
\special{pa 2437 490}%
\special{pa 2457 476}%
\special{pa 2477 490}%
\special{pa 2457 423}%
\special{fp}%
%
\special{pn 20}%
\special{pa 1880 529}%
\special{pa 1880 423}%
\special{fp}%
\special{sh 1}%
\special{pa 1880 423}%
\special{pa 1860 490}%
\special{pa 1880 476}%
\special{pa 1900 490}%
\special{pa 1880 423}%
\special{fp}%
%
\special{pn 20}%
\special{pa 933 1180}%
\special{pa 955 1159}%
\special{pa 978 1139}%
\special{pa 1000 1118}%
\special{pa 1023 1097}%
\special{pa 1046 1077}%
\special{pa 1069 1057}%
\special{pa 1093 1037}%
\special{pa 1117 1017}%
\special{pa 1142 998}%
\special{pa 1167 978}%
\special{pa 1193 959}%
\special{pa 1219 941}%
\special{pa 1246 923}%
\special{pa 1274 905}%
\special{pa 1303 888}%
\special{pa 1333 871}%
\special{pa 1363 855}%
\special{pa 1395 839}%
\special{pa 1428 824}%
\special{pa 1462 810}%
\special{pa 1496 796}%
\special{pa 1531 782}%
\special{pa 1566 768}%
\special{pa 1600 754}%
\special{pa 1634 739}%
\special{pa 1667 725}%
\special{pa 1699 709}%
\special{pa 1729 693}%
\special{pa 1757 675}%
\special{pa 1783 657}%
\special{pa 1806 637}%
\special{pa 1826 615}%
\special{pa 1843 592}%
\special{pa 1856 567}%
\special{pa 1866 540}%
\special{pa 1871 511}%
\special{pa 1872 479}%
\special{pa 1871 452}%
\special{sp}%
%
\special{pn 20}%
\special{pa 1871 1190}%
\special{pa 1850 1165}%
\special{pa 1829 1140}%
\special{pa 1807 1115}%
\special{pa 1785 1091}%
\special{pa 1763 1068}%
\special{pa 1740 1046}%
\special{pa 1716 1025}%
\special{pa 1691 1005}%
\special{pa 1666 987}%
\special{pa 1639 971}%
\special{pa 1611 956}%
\special{pa 1582 943}%
\special{pa 1552 932}%
\special{pa 1522 921}%
\special{pa 1491 911}%
\special{pa 1460 901}%
\special{pa 1450 898}%
\special{sp}%
%
\special{pn 20}%
\special{pa 959 500}%
\special{pa 933 432}%
\special{fp}%
\special{sh 1}%
\special{pa 933 432}%
\special{pa 938 501}%
\special{pa 952 482}%
\special{pa 975 487}%
\special{pa 933 432}%
\special{fp}%
%
\special{pn 20}%
\special{pa 959 491}%
\special{pa 971 522}%
\special{pa 982 552}%
\special{pa 995 582}%
\special{pa 1009 610}%
\special{pa 1025 638}%
\special{pa 1042 663}%
\special{pa 1062 688}%
\special{pa 1084 710}%
\special{pa 1108 731}%
\special{pa 1133 751}%
\special{pa 1159 770}%
\special{pa 1186 789}%
\special{pa 1214 806}%
\special{pa 1242 824}%
\special{pa 1252 830}%
\special{sp}%
\put(17.5100,-3.5500){\makebox(0,0)[lb]{$t_i'$}}%
\put(17.8000,-14.7000){\makebox(0,0)[lb]{$t_j$}}%
\put(23.9600,-14.7100){\makebox(0,0)[lb]{$t_k$}}%
\put(4.3000,-8.8000){\makebox(0,0)[lb]{$\b_i$}}%
\put(8.5600,-3.4500){\makebox(0,0)[lb]{$t_j'$}}%
\put(23.8800,-3.4500){\makebox(0,0)[lb]{$t_k'$}}%
\put(8.4700,-14.8100){\makebox(0,0)[lb]{$t_i$}}%
%
\special{pn 13}%
\special{pa 3760 805}%
\special{pa 6068 795}%
\special{fp}%
%
\special{pn 20}%
\special{ar 4550 805 430 430  1.5946014 6.2831853}%
\special{ar 4550 805 430 430  0.0000000 1.5707963}%
%
\special{pn 20}%
\special{sh 0.600}%
\special{ar 4980 810 36 36  0.0000000 6.2831853}%
%
\special{pn 20}%
\special{sh 0.600}%
\special{ar 4120 800 35 35  0.0000000 6.2831853}%
%
\special{pn 20}%
\special{sh 0.600}%
\special{ar 5710 790 36 36  0.0000000 6.2831853}%
%
\special{pn 20}%
\special{pa 4960 935}%
\special{pa 4970 885}%
\special{fp}%
\special{sh 1}%
\special{pa 4970 885}%
\special{pa 4937 946}%
\special{pa 4960 937}%
\special{pa 4977 954}%
\special{pa 4970 885}%
\special{fp}%
%
\special{pn 20}%
\special{pa 4590 380}%
\special{pa 4530 380}%
\special{fp}%
\special{sh 1}%
\special{pa 4530 380}%
\special{pa 4597 400}%
\special{pa 4583 380}%
\special{pa 4597 360}%
\special{pa 4530 380}%
\special{fp}%
%
\special{pn 20}%
\special{pa 4530 1235}%
\special{pa 4600 1235}%
\special{fp}%
\special{sh 1}%
\special{pa 4600 1235}%
\special{pa 4533 1215}%
\special{pa 4547 1235}%
\special{pa 4533 1255}%
\special{pa 4600 1235}%
\special{fp}%
\put(56.6000,-10.7000){\makebox(0,0)[lb]{$t_k$}}%
\put(50.2000,-10.8000){\makebox(0,0)[lb]{$t_j$}}%
\put(39.5000,-10.4500){\makebox(0,0)[lb]{$t_i$}}%
\put(60.7000,-4.3000){\makebox(0,0)[lb]{$\hat{\C}$}}%
%
\special{pn 13}%
\special{pa 3550 180}%
\special{pa 6470 180}%
\special{pa 6470 1460}%
\special{pa 3550 1460}%
\special{pa 3550 180}%
\special{fp}%
\put(61.5000,-8.7000){\makebox(0,0)[lb]{$\hat{\bR}$}}%
\put(1.2000,-4.6000){\makebox(0,0)[lb]{$T$}}%
\put(1.3000,-12.4000){\makebox(0,0)[lb]{$T$}}%
%
\special{pn 20}%
\special{pa 4130 680}%
\special{pa 4120 730}%
\special{fp}%
\special{sh 1}%
\special{pa 4120 730}%
\special{pa 4153 669}%
\special{pa 4130 678}%
\special{pa 4113 661}%
\special{pa 4120 730}%
\special{fp}%
\end{picture}%
\end{center}
\caption{Basic braid $\b_i$ in $T$ and the corresponding 
movement of $t$ in $\hat{\C}$}
\label{fig:braid}
\end{figure}

 
Let us describe the fundamental group $\pi_1(T,t)$ in terms of 
a braid group \cite{Birman}. 
We take a base point $t = (t_1,t_2,t_3) \in T$ in such a manner 
that the three points lie on the real line in an increasing order 
$t_1 < t_2 < t_3$.  
To treat them symmetrically, we denote them by $t_i$, $t_j$, $t_k$ 
for a cyclic permutation $(i,j,k)$ of $(1,2,3)$ and think of 
them as cyclically ordered three points on the equator 
$\hat{\bR} = \bR \cup \{\infty\}$ of the Riemann sphere 
$\hat{\C} = \C \cup \{\infty\}$. 
Let $\b_i$ be a braid on three strings as in 
Figure \ref{fig:braid} (left) along which $t_i$ and $t_j$ make a 
half-turn, with $t_i$ moving in the southern hemisphere and 
$t_j$ in the northern hemisphere, while $t_k$ is kept fixed as in 
Figure \ref{fig:braid} (right). 
Then the braid group on three strings is the group generated by 
$\b_i$, $\b_j$ and $\b_k$ and the pure braid group $P_3$ is the 
normal subgroup of $B_3$ generated by the squares 
$\b_i^2$, $\b_j^2$ and $\b_k^2$, 
\[
P_3 = \la \b_i^2, \b_j^2, \b_k^2 \ra \triangleleft 
B_3 = \la \b_i, \b_j, \b_k \ra. 
\]
The generators of $B_3$ satisfy relations 
$\b_i\b_j\b_i = \b_j\b_i\b_j$ and $\b_k = \b_i \b_j \b_i^{-1}$, 
and so $B_3$ is generated by $\b_i$ and $\b_j$ only. 
The fundamental group $\pi_1(T,t)$ is identified with 
the pure braid group $P_3$. 
\par 
The reduction map (\ref{eqn:reduction}) induces a group 
homomorphism $r_* : \pi_1(T,t) \to \pi_1(X,x)$. 
It is easy to see that the loops $\ell_0$ and $\ell_1$ in 
Figure \ref{fig:loop} are the $r_*$-images of $\b_1^2$ and 
$\b_2^2$ respectively, so that the Pochhammer loop $\wp$ in $X$ 
is the $r_*$-image of the pure braid 
\begin{equation} \label{eqn:commutator1} 
[\b_1^2, \b_2^{-2}] = \b_1^2 \b_2^{-2} \b_1^{-2} \b_2^2. 
\end{equation}
Thus we will be concerned with the Poincar\'e section 
(\ref{eqn:PSt}) along this particular braid. 
\par 
The symmetric group $S_3$ acts on $T$ by permuting the entries 
of $t = (t_1,t_2,t_3)$ and the quotient space 
$T/S_3$ is the configuration space of mutually distinct, 
unordered, three points in $\C$. 
The fundamental group $\pi(T/S_3, s)$ with base point 
$s = \{t_1,t_2,t_3\}$ is identified with the ordinary braid group 
$B_3$ and there exists a short exact sequence of groups 
\[
\begin{CD}
1 @>>> \pi_1(T,t) @>>> \pi_1(T/S_3,s) @>>> S_3 @>>> 1 \\
@.       @|                    @|           @|    @.  \\
1 @>>> P_3        @>>> B_3              @>>> S_3 @>>> 1.
\end{CD}
\]
Then the Poincar\'e section (\ref{eqn:PSt}) naturally lifts to 
a collection of isomorphisms 
\[
\b_* : \M_t(\k) \to \M_{\tau(t)}(\tau(\k)), \qquad (\b \in B_3) 
\]
which may be called the {\sl half-Poincar\'e section} of 
$\PVI(\k)$, where $\tau \in S_3$ denotes the permutation 
corresponding to $\b \in B_3$. 
Note that $\tau \in S_3$ acts on $\k \in \K$ too by permuting 
the entries of $(\k_1,\k_2,\k_3)$ in the same manner as it 
does on $t = (t_1,t_2,t_3)$, since $\k_i$ is loaded on $t_i$. 
Now the permutation corresponding to the basic braid $\b_i$ is 
the substitution $\tau_i = (i,j)$ that exchanges $t_i$ and 
$t_j$ while keeping $t_k$ fixed. 
Thus there are three basic half-Poincar\'e maps: 
\begin{equation} \label{eqn:PSi}
\b_{i*} : \M_t(\k) \to \M_{\tau_i(t)}(\tau_i(\k)), 
\qquad (i = 1,2,3).  
\end{equation}
\section{Riemann-Hilbert Correspondence} \label{sec:RHC} 
It is rather hopeless to deal with the Painlev\'e flow directly, 
since it is a highly transcendental dynamical system 
on the moduli space of stable parabolic connections. 
But it can be recast into a more tractable dynamical system, 
called an isomonodromic flow, on a moduli space of monodromy 
representations via a Riemann-Hilbert correspondence. 
We review the construction of such a Riemann-Hilbert 
correspondence in the sequel. 
\par
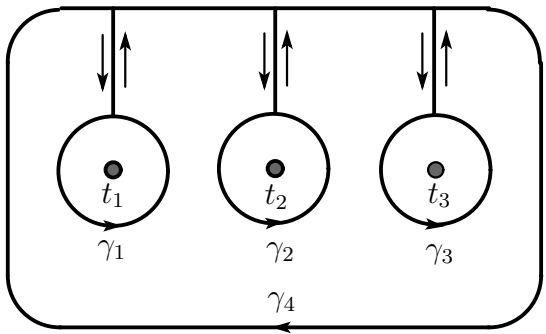
\begin{figure}[t]
\begin{center}
\unitlength 0.1in
\begin{picture}(27.93,16.78)(2.10,-18.78)
%
\special{pn 20}%
\special{ar 483 480 273 273  3.1165979 4.6873942}%
%
\special{pn 20}%
\special{pa 210 480}%
\special{pa 210 1600}%
\special{fp}%
%
\special{pn 20}%
\special{ar 483 1600 273 273  1.5458015 3.1415927}%
%
\special{pn 20}%
\special{ar 1603 1040 287 287  0.0000000 6.2831853}%
%
\special{pn 20}%
\special{ar 2450 1047 287 287  0.0000000 6.2831853}%
%
\special{pn 20}%
\special{pa 490 200}%
\special{pa 2730 200}%
\special{fp}%
%
\special{pn 20}%
\special{pa 483 1873}%
\special{pa 2730 1873}%
\special{fp}%
%
\special{pn 20}%
\special{ar 2730 1600 273 273  6.2831853 6.2831853}%
\special{ar 2730 1600 273 273  0.0000000 1.5964317}%
%
\special{pn 20}%
\special{pa 763 200}%
\special{pa 763 760}%
\special{fp}%
%
\special{pn 20}%
\special{pa 1610 200}%
\special{pa 1610 753}%
\special{fp}%
%
\special{pn 20}%
\special{pa 2443 200}%
\special{pa 2443 753}%
\special{fp}%
%
\special{pn 20}%
\special{pa 3003 487}%
\special{pa 3003 1607}%
\special{fp}%
%
\special{pn 20}%
\special{ar 2723 480 274 274  4.7123890 6.2831853}%
\special{ar 2723 480 274 274  0.0000000 0.0256354}%
%
\special{pn 20}%
\special{ar 763 1054 287 287  0.0000000 6.2831853}%
%
\special{pn 13}%
\special{sh 0.600}%
\special{ar 2450 1047 41 41  0.0000000 6.2831853}%
%
\special{pn 20}%
\special{sh 0.600}%
\special{ar 1610 1040 41 41  0.0000000 6.2831853}%
%
\special{pn 20}%
\special{sh 0.600}%
\special{ar 763 1047 41 41  0.0000000 6.2831853}%
%
\special{pn 13}%
\special{pa 707 340}%
\special{pa 707 620}%
\special{fp}%
\special{sh 1}%
\special{pa 707 620}%
\special{pa 727 553}%
\special{pa 707 567}%
\special{pa 687 553}%
\special{pa 707 620}%
\special{fp}%
%
\special{pn 13}%
\special{pa 1554 333}%
\special{pa 1554 613}%
\special{fp}%
\special{sh 1}%
\special{pa 1554 613}%
\special{pa 1574 546}%
\special{pa 1554 560}%
\special{pa 1534 546}%
\special{pa 1554 613}%
\special{fp}%
%
\special{pn 13}%
\special{pa 2380 333}%
\special{pa 2380 613}%
\special{fp}%
\special{sh 1}%
\special{pa 2380 613}%
\special{pa 2400 546}%
\special{pa 2380 560}%
\special{pa 2360 546}%
\special{pa 2380 613}%
\special{fp}%
%
\special{pn 13}%
\special{pa 826 613}%
\special{pa 826 340}%
\special{fp}%
\special{sh 1}%
\special{pa 826 340}%
\special{pa 806 407}%
\special{pa 826 393}%
\special{pa 846 407}%
\special{pa 826 340}%
\special{fp}%
%
\special{pn 13}%
\special{pa 1673 606}%
\special{pa 1673 333}%
\special{fp}%
\special{sh 1}%
\special{pa 1673 333}%
\special{pa 1653 400}%
\special{pa 1673 386}%
\special{pa 1693 400}%
\special{pa 1673 333}%
\special{fp}%
%
\special{pn 13}%
\special{pa 2506 606}%
\special{pa 2506 333}%
\special{fp}%
\special{sh 1}%
\special{pa 2506 333}%
\special{pa 2486 400}%
\special{pa 2506 386}%
\special{pa 2526 400}%
\special{pa 2506 333}%
\special{fp}%
%
\special{pn 13}%
\special{pa 728 1341}%
\special{pa 791 1341}%
\special{fp}%
\special{sh 1}%
\special{pa 791 1341}%
\special{pa 724 1321}%
\special{pa 738 1341}%
\special{pa 724 1361}%
\special{pa 791 1341}%
\special{fp}%
%
\special{pn 13}%
\special{pa 1568 1327}%
\special{pa 1638 1320}%
\special{fp}%
\special{sh 1}%
\special{pa 1638 1320}%
\special{pa 1570 1307}%
\special{pa 1585 1325}%
\special{pa 1574 1347}%
\special{pa 1638 1320}%
\special{fp}%
%
\special{pn 13}%
\special{pa 2422 1334}%
\special{pa 2478 1341}%
\special{fp}%
\special{sh 1}%
\special{pa 2478 1341}%
\special{pa 2414 1313}%
\special{pa 2425 1334}%
\special{pa 2409 1353}%
\special{pa 2478 1341}%
\special{fp}%
%
\special{pn 13}%
\special{pa 1673 1873}%
\special{pa 1610 1873}%
\special{fp}%
\special{sh 1}%
\special{pa 1610 1873}%
\special{pa 1677 1893}%
\special{pa 1663 1873}%
\special{pa 1677 1853}%
\special{pa 1610 1873}%
\special{fp}%
\put(7.5600,-14.6700){\makebox(0,0){$\gamma_1$}}%
\put(16.3800,-14.8100){\makebox(0,0){$\gamma_2$}}%
\put(24.7100,-14.8800){\makebox(0,0){$\gamma_3$}}%
\put(16.4500,-17.4000){\makebox(0,0){$\gamma_4$}}%
\put(7.6300,-11.8000){\makebox(0,0){$t_1$}}%
\put(16.1700,-11.8700){\makebox(0,0){$t_2$}}%
\put(24.7100,-11.8000){\makebox(0,0){$t_3$}}%
\end{picture}%
\end{center}
\caption{Four loops in $\P^1-D_t$; 
the fourth point $t_4$ is outside $\ga_4$, invisible.} 
\label{fig:loop2}
\end{figure}

 
Let $A := \C^4$ be the complex $4$-space with 
coordinates $a = (a_1,a_2,a_3,a_4)$, called the space of 
local monodromy data. 
Given $(t,a) \in T \times A$, let $\R_t(a)$ be the moduli 
space of Jordan equivalence classes of representations 
$\rho : \pi_1(\P^1-D_t,*) \to SL_2(\C)$
such that $\Tr\, \rho(\ga_i) = a_i$ for $i \in \{1,2,3,4\}$, 
where the divisor $D_t = t_1+t_2+t_3+t_4$ is identified with 
the point set $\{t_1,t_2,t_3,t_4\}$ and $\ga_i$ is a loop 
surrounding $t_i$ as in Figure \ref{fig:loop2}. 
Any stable parabolic connection $Q = (E,\nabla,\psi,l) 
\in \M_t(\k)$, when restricted to $\P^1-D_t$, induces a flat 
connection 
\[
\nabla|_{\P^1-D_t} : E|_{\P^1-D_t} \to (E|_{\P^1-D_t}) 
\otimes \Omega_{\P^1-D_t}^1,
\]
and one can speak of the Jordan equivalence class $\rho$ of 
its monodromy representations. 
Then the Riemann-Hilbert correspondence at $t \in T$ is 
defined by  
\begin{equation} \label{eqn:RHtk}
\RH_{t,\k} : \M_t(\k) \to \R_t(a), \quad Q \mapsto \rho, 
\end{equation} 
where in view of the Riemann scheme in Table \ref{tab:riemann}, 
the local monodromy data $a \in A$ is given by 
\begin{equation} \label{eqn:a}
a_i = \left\{
\begin{array}{ll}
{\-}2 \cos \pi \k_i \qquad & (i = 1,2,3), \\[2mm] 
-2 \cos \pi \k_4 \qquad & (i = 4). 
\end{array}
\right.
\end{equation}
\par 
As a relative setting over $T$, let $\pi_a : \R(a) \to T$ be 
the family of moduli spaces of monodromy representations with 
fiber $\R_t(a)$ over $t \in T$. 
Then the relative version of Riemann-Hilbert correspondence 
is formulated to be the commutative diagram 
\begin{equation} \label{eqn:RHk}
\begin{CD} 
\M(\k) @> \RH_{\k} >> \R(a)    \\
@V \pi_{\k} VV     @VV \pi_a V \\
T @=  T, 
\end{CD}
\end{equation}
whose fiber over $t \in T$ is given by (\ref{eqn:RHtk}). 
Then we have the following theorem \cite{IIS1,IIS2}. 
\begin{theorem} \label{thm:RH} 
If $\k \in \K-\Wall$, then $\R(a)$ as well as each fiber 
$\R_t(a)$ is smooth and the Riemann-Hilbert correspondence 
$\RH_{\k}$ in $(\ref{eqn:RHk})$ is a biholomorphism. 
\end{theorem}
\begin{remark} \label{rem:RH} 
If $\k \in \Wall$, then $\R_t(a)$ is not a smooth surface 
but a surface with Klein singularities and (\ref{eqn:RHtk}) 
yields an analytic minimal resolution of singularities, 
so that (\ref{eqn:RHk}) gives a family of resolutions of 
singularities \cite{IIS1}. 
As is mentioned in Remark \ref{rem:main2}, this fact makes 
the treatment of the nongeneric case more involved and 
we leave this case in another occasion. 
\end{remark}
\section{Cubic Surface and the 27 Lines} \label{sec:cubic}
The moduli space $\R_t(a)$ of monodromy representations 
is isomorphic to an affine cubic surface $\Sol(\th)$ and 
the braid group action on $\R_t(a)$ can be made explicit 
in terms of $\Sol(\th)$. 
Let us recall this construction \cite{IIS1}. 
Given $\th = (\th_1,\th_2,\th_3,\th_4) \in \Th := \C^4_{\th}$, 
consider an affine cubic surface 
\[
\Sol(\th) = \{\, x = (x_1,x_2,x_3) \in \C^3_x \,:\, 
f(x,\th) = 0 \, \}, 
\]
where the cubic polynomial $f(x,\th)$ of $x$ with parameter 
$\th$ is given by 
\[
f(x,\th) = x_1x_2x_3 + x_1^2 + x_2^2 + x_3^2 
- \th_1 x_1 - \th_2 x_2 - \th_3 x_3 + \th_4. 
\]
Then there exists an isomorphism of affine algebraic 
surfaces, $\R_t(a) \to \Sol(\th)$, $\rho \mapsto x$, 
where 
\[
x_i = \Tr\,\rho(\ga_j\ga_k), \qquad 
\text{for} \quad \{i,j,k\} = \{1,2,3\}, 
\]
together with a correspondence of parameters, 
$A \to \Th$, $a \mapsto \th$, given by 
\begin{equation} \label{eqn:th}
\th_i = \left\{
\begin{array}{ll}
a_i a_4 + a_j a_k \qquad & (\{i,j,k\}=\{1,2,3\}), \\[2mm]
a_1 a_2 a_3 a_4 + a_1^2 + a_2^2 + a_3^2 + a_4^2 - 4 
\qquad & (i = 4). 
\end{array}
\right.
\end{equation}
\par 
With this identification, the Riemann-Hilbert correspondence 
(\ref{eqn:RHtk}) is reformulated as a map
\begin{equation} \label{eqn:RHtk2}
\RH_t(\k) : \M_t(\k) \to \Sol(\th), 
\qquad \text{with} \quad \th = \rh(\k), 
\end{equation}
where $\rh : \K \to \Th$ is the composition of two maps 
$\K \to A$ and $A \to \Th$ defined by (\ref{eqn:a}) and 
(\ref{eqn:th}), which we call the Riemann-Hilbert correspondence 
in the parameter level. 
Through the reformulated Riemann-Hilbert correspondence 
(\ref{eqn:RHtk2}), the $i$-th basic half-Poincar\'e map $\b_{i*}$ 
in (\ref{eqn:PSi}) is conjugated to a map 
$g_i : \Sol(\th) \to \Sol(\th')$, $(x,\th) \mapsto (x',\th')$, 
which is explicitly represented as 
\begin{equation} \label{eqn:gi}
g_i \quad : \quad 
(x_i',x_j',x_k',\th_i',\th_j',\th_k',\th_4') = 
(\th_j-x_j-x_kx_i, x_i, x_k, \th_j, \th_i, \th_k, \th_4). 
\end{equation}
A derivation of this formula can be found in \cite{Iwasaki4} 
(see also \cite{Boalch2,DM,Goldman,Guzzetti,Iwasaki3,Jimbo}). 
By Theorem \ref{thm:RH} the map (\ref{eqn:RHtk2}) is an 
isomorphism and hence (\ref{eqn:gi}) is a strict conjugacy 
of (\ref{eqn:PSi}). 
We can easily check the relations $g_ig_jg_i = g_jg_ig_j$ and 
$g_k = g_i g_j g_i^{-1}$ which are parallel to those for 
$\b_i$, $\b_j$, $\b_k$. 
\par
\begin{figure}[t] 
\begin{center}
\unitlength 0.1in
\begin{picture}(20.68,21.44)(14.16,-27.54)
%
\special{pn 13}%
\special{ar 2450 1720 1034 1034  4.9299431 6.2831853}%
\special{ar 2450 1720 1034 1034  0.0000000 4.4841152}%
%
\special{pn 20}%
\special{pa 1834 2000}%
\special{pa 3104 2000}%
\special{fp}%
%
\special{pn 20}%
\special{pa 2644 1020}%
\special{pa 1934 2210}%
\special{fp}%
%
\special{pn 20}%
\special{pa 2254 1020}%
\special{pa 2964 2200}%
\special{fp}%
\put(22.9000,-7.8000){\makebox(0,0)[lb]{$\ol{\Sol}(\th)$}}%
\put(29.6000,-14.5000){\makebox(0,0)[lb]{$\Sol(\th)$}}%
\put(23.7000,-22.1000){\makebox(0,0)[lb]{$L_i$}}%
\put(20.3000,-17.0000){\makebox(0,0)[lb]{$L_j$}}%
\put(27.1000,-16.9000){\makebox(0,0)[lb]{$L_k$}}%
\put(30.2000,-21.8000){\makebox(0,0)[lb]{$p_j$}}%
\put(23.9000,-11.9000){\makebox(0,0)[lb]{$p_i$}}%
\put(17.9000,-21.7000){\makebox(0,0)[lb]{$p_k$}}%
%
\special{pn 20}%
\special{sh 0.600}%
\special{ar 2060 2000 40 40  0.0000000 6.2831853}%
%
\special{pn 20}%
\special{sh 0.600}%
\special{ar 2450 1340 40 40  0.0000000 6.2831853}%
%
\special{pn 20}%
\special{sh 0.600}%
\special{ar 2840 2000 40 40  0.0000000 6.2831853}%
\end{picture}%
\end{center}
\caption{Tritangent lines at infinity on $\ol{\Sol}(\th)$} 
\label{fig:cubic1} 
\end{figure}
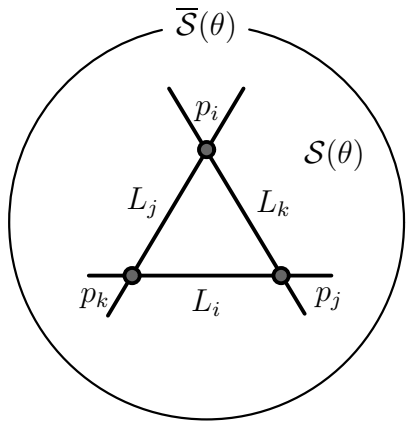

 
To utilize standard techniques from algebraic geometry and 
complex geometry, we need to compactify the affine cubic 
surface $\Sol(\th)$ by a standard embedding 
\[
\Sol(\th) \hookrightarrow \ol{\Sol}(\th) \subset \P^3, \qquad 
x = (x_1,x_2,x_3) \mapsto [1:x_1:x_2:x_3], 
\]
where the compactified surface $\ol{\Sol}(\th)$ is defined 
by $\ol{\Sol}(\th) = \{\,X \in \P^3\,:\, F(X,\th)=0 \,\}$ with 
\[
F(X,\th) = X_1X_2X_3+X_0(X_1^2+X_2^2+X_3^2)-
X_0^2(\th_1X_1+\th_2X_2+\th_3X_3)+\th_4X_0^3. 
\]
It is obtained from the affine surface $\Sol(\th)$ by adding 
three lines at infinity, 
\begin{equation} \label{eqn:lines}
L_i = \{\, X \in \P^3 \,:\, X_0 = X_i = 0 \,\} \qquad 
(i = 1,2,3). 
\end{equation}
Here and hereafter the homogeneous coordinates 
$X = [X_0:X_1:X_2:X_3]$ of $\P^3$ should not be confused with 
the domain $X$ in (\ref{eqn:X}). 
The union $L = L_1 \cup L_2 \cup L_3$ is called the 
tritangent lines at infinity and the intersection point of 
$L_j$ and $L_k$ is denoted by $p_i$ (see Figure \ref{fig:cubic1}). 
Note that 
\[
p_1 = [0:1:0:0], \qquad p_2 = [0:0:1:0], \qquad 
p_3 = [0:0:0:1]. 
\]
For $i \in \{1,2,3\}$, put 
$U_i = \{\, X \in \P^3 \,:\, X_i \neq 0 \,\}$ and 
take inhomogeneous coordinates of $\P^3$; 
\begin{equation} \label{eqn:coord}
\begin{array}{rclcll}
u &=& (u_0,u_j,u_k) &=& (X_0/X_i, X_j/X_i, X_k/X_i) \qquad & 
\text{on} \quad U_i, \\[1mm]
v &=& (v_0,v_i,v_k) &=& (X_0/X_j, X_i/X_j, X_k/X_j) \qquad & 
\text{on} \quad U_j, \\[1mm]
w &=& (w_0,w_i,w_j) &=& (X_0/X_k, X_i/X_k, X_j/X_k) \qquad & 
\text{on} \quad U_k, 
\end{array}
\end{equation}
where $\{i,j,k\} = \{1,2,3\}$. 
In terms of these coordinates we shall find local coordinates 
and local equations of $\ol{\Sol}(\th)$ around $L$. 
Since $L \subset U_1 \cup U_2 \cup U_3$, we can divide $L$ 
into components $L \cap U_i$, $i = 1,2,3$, and make a further 
decomposition $L \cap U_i = \{p_i\} \cup (L_j-\{p_i,p_k\}) 
\cup (L_k-\{p_i,p_j\})$ into a total of nine pieces. 
Then a careful inspection of equation $F(X,\th) = 0$ implies 
that around those pieces we can take local coordinates and 
local equations as in Table \ref{tab:local}, where 
$O_m(u_j, u_k) = O((|u_j|+|u_k|)^m)$ denotes a small term of 
order $m$ as $(u_j, u_k) \to (0,0)$. 
\begin{table}[t] 
\begin{center}
\begin{tabular}{|c|c|l|}
\hline
\vspace{-0.3cm} & &  \\
coordinates & valid around & local equation \\[2mm]
\hline
\hline
\vspace{-0.2cm} & &  \\
$(u_j,u_k)$ & $p_i$ & 
$u_0 = -(u_j u_k) \{1-(u_j^2+\th_i u_j u_k+u_k^2)+O_3(u_j,u_k)\}$ \\[2mm]
\hline
\vspace{-0.2cm} & &  \\
$(u_0,u_k)$ & $L_j-\{ p_i,p_k \}$ & 
$u_j = -(u_k+1 / u_k)u_0+(\th_k+\th_i / u_k) u_0^2 + O(u_0^3)$ \\[2mm]
\hline
\vspace{-0.2cm} & &  \\
$(u_0,u_j)$ & $L_k-\{ p_i,p_j \} $ & 
$u_k = -(u_j+1 / u_j)u_0+(\th_j+\th_i / u_j) u_0^2 + O(u_0^3)$ \\[2mm]
\hline
\hline
\vspace{-0.2cm} & &  \\
$(v_i,v_k)$ & $p_j$ & 
$v_0 = -(v_i v_k) \{1-(v_i^2+\th_j v_i v_k+v_k^2)+O_3(v_i,v_k)\}$ \\[2mm]
\hline
\vspace{-0.2cm} & &  \\
$(v_0,v_i)$ & $L_k - \{ p_i,p_j \}$ & 
$v_k = -(v_i+1 / v_i)v_0+(\th_i+\th_j / v_i) v_0^2 + O(v_0^3)$ \\[2mm]
\hline
\vspace{-0.2cm} & &  \\
$(v_0,v_k)$ & $L_i - \{ p_j,p_k \}$ & 
$v_i = -(v_k+1 / v_k)v_0+(\th_k+\th_j / v_k) v_0^2 + O(v_0^3)$ \\[2mm]
\hline
\hline
\vspace{-0.2cm} & &  \\
$(w_i,w_j)$ & $p_k$ & 
$w_0 = -(w_i w_j) \{1-(w_i^2+\th_k w_i w_j+w_j^2)+O_3(w_i,w_j)\}$ \\[2mm]
\hline
\vspace{-0.2cm} & &  \\
$(w_0,w_j)$ & $L_i -\{ p_j,p_k\}$ & 
$w_i = -(w_j+1 / w_j)w_0+(\th_j+\th_k / w_j) w_0^2 + O(w_0^3)$ \\[2mm]
\hline
\vspace{-0.2cm} & &  \\
$(w_0,w_i)$ & $L_j - \{ p_i,p_k \}$ & 
$w_j = -(w_i+1 / w_i)w_0+(\th_i+\th_k / w_i) w_0^2 + O(w_0^3)$ \\[2mm]
\hline
\end{tabular}
\end{center}
\caption{Local coordinates and local equations of $\ol{\Sol}(\th)$}
\label{tab:local}
\end{table}
\begin{lemma} \label{lem:smooth} 
As to the smoothness of the surface $\ol{\Sol}(\th)$, 
the following hold. 
\begin{enumerate} 
\item For any $\th \in \Th$, the surface $\ol{\Sol}(\th)$ is 
smooth in a neighborhood of $L$.  
\item If $\th = \rh(\k)$ with $\k \in \K$, the surface 
$\ol{\Sol}(\th)$ is smooth everywhere if and only if 
$\k \in \K-\Wall$. 
\end{enumerate} 
\end{lemma}
{\it Proof}. 
First we show assertion (1). 
In terms of inhomogeneous coordinates $u$ in (\ref{eqn:coord}), 
we have 
\[
\ol{\Sol}(\th) \cap U_i \cong 
\{\, u = (u_0,u_j,u_k) \in \C^3 \,:\, f_i(u,\th) = 0 \,\},
\]
where the defining equation $f_i(u,\th)$ is given by 
\[
f_i(u,\th) = u_j u_k +u_0 (1 + u_j^2 + u_k^2) - u_0^2 
(\th_i + \th_j u_j + \th_k u_k) + \th_4 u_0^3. 
\]
The partial derivatives of $f_i = f_i(u,\th)$ with respect 
to $u = (u_0, u_j, u_k)$ are calculated as 
\begin{eqnarray*}
\frac{\partial f_i}{\partial u_0} 
&=& (1 + u_j^2 + u_k^2) 
- 2 u_0 (\th_i + \th_j u_j + \th_k u_k) + 3 \th_4 u_0^2 \\[2mm]
\frac{\partial f_i}{\partial u_j} 
&=& u_k + 2 u_0 u_j - \th_j u_0^2 \\[2mm]
\frac{\partial f_i}{\partial u_k} 
&=& u_j + 2 u_0 u_k - \th_k u_0^2 .
\end{eqnarray*}
Restricted to $L \cap U_i = (L_j \cap U_i) \cup (L_k \cap U_i)$, 
these derivatives become 
\[
\begin{array}{rclrclrcll}
\dfrac{\partial f_i}{\partial u_0} &=& 1 + u_k^2, \qquad & 
\dfrac{\partial f_i}{\partial u_j} &=& u_k,       \qquad & 
\dfrac{\partial f_i}{\partial u_k} &=& 0,         \qquad & 
\text{on} \quad L_j \cap U_i, \\[4mm]
\dfrac{\partial f_i}{\partial u_0} &=& 1 + u_j^2, \qquad & 
\dfrac{\partial f_i}{\partial u_j} &=& 0,         \qquad & 
\dfrac{\partial f_i}{\partial u_k} &=& u_j,       \qquad & 
\text{on} \quad L_k \cap U_i. 
\end{array}
\]
Hence the exterior derivative $d_u f_i$ does not vanish on 
$L \cap U_i$, and the implicit function theorem implies that 
$\ol{\Sol}(\th)$ is smooth in a neighborhood of $L$. 
This proves assertion (1). 
In order to show assertion (2) we recall that the affine 
surface $\Sol(\th)$ is smooth if and only if $\th = \rh(\k)$ 
with $\k \in \K-\Wall$ (see \cite{IIS1}). 
Then assertion (2) readily follows from assertion (1). 
\hfill $\Box$ 
\par\medskip 
Now let us review some basic facts about smooth cubic surfaces 
in $\P^3$ (see e.g. \cite{GH}).
It is well known that every smooth cubic surface $S$ 
in $\P^3$ can be obtained by blowing up $\P^2$ at six points 
$P_1, \dots, P_6$, no three colinear and not all six on a conic, 
and embedding the blow-up surface into $\P^3$ by the proper 
transform of the linear system of cubics passing through the 
six points $P_1, \dots, P_6$. 
It is also well known that there are exactly $27$ 
lines on the smooth cubic surface $S$, each of which has 
self-intersection number $-1$. 
Explicitly, they are given by
\[
E_a \quad (a = 1,\dots,6); \qquad F_{ab} \quad 
(1 \le a < b \le 6); \qquad G_a \quad (a = 1,\dots,6), 
\]
\begin{enumerate}
\item $E_a$ is the exceptional curve over the point $P_a$, 
\item $F_{ab}$ is the strict transform of the line in 
$\P^2$ through the two points $P_a$ and $P_b$, 
\item $G_a$ is the strict transform of the conic in $\P^2$ 
through the five points $P_1,\dots,\hat{P}_a,\dots,P_6$. 
\end{enumerate}
Here the index $a$ should not be confused with the 
local monodromy data $a \in A$.  
All the intersection relations among the $27$ lines 
with {\sl nonzero} intersection numbers are listed as 
\[
\begin{array}{rl}
(E_a, E_a) = (G_a, G_a) = (F_{ab}, F_{ab}) = -1 \qquad & 
(\forall \, a,b), \\[2mm]
(E_a, F_{bc}) = (G_a, F_{bc}) = {\-} 1  \qquad & 
(a \in \{b,c\}), \\[2mm] 
(E_a, G_b) = {\-} 1 \qquad & 
(a \neq b),  \\[2mm]
(F_{ab}, F_{cd}) = {\-} 1 \qquad & 
(\{a,b\} \cap \{c,d\} = \emptyset).
\end{array}
\]
\par
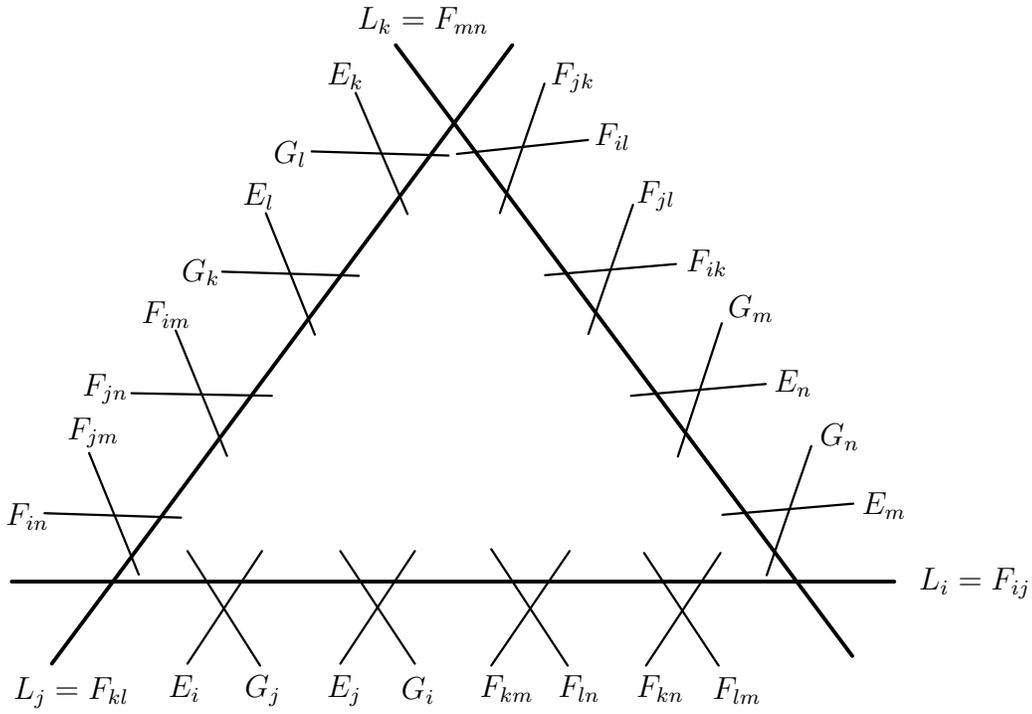
\begin{figure}[t] 
\begin{center}
\unitlength 0.1in
\begin{picture}(47.80,35.30)(11.40,-35.70)
\put(59.2000,-31.7000){\makebox(0,0)[lb]{$L_i=F_{ij}$}}%
\put(11.8000,-37.4000){\makebox(0,0)[lb]{$L_j=F_{kl}$}}%
\put(29.8000,-2.1000){\makebox(0,0)[lb]{$L_k=F_{mn}$}}%
%
\special{pn 20}%
\special{pa 1170 3080}%
\special{pa 5790 3080}%
\special{fp}%
%
\special{pn 20}%
\special{pa 3790 270}%
\special{pa 1380 3510}%
\special{fp}%
\special{pa 3180 270}%
\special{pa 5570 3470}%
\special{fp}%
%
\special{pn 13}%
\special{pa 2090 2920}%
\special{pa 2470 3520}%
\special{fp}%
\special{pa 2480 2920}%
\special{pa 2090 3510}%
\special{fp}%
\special{pa 2890 2920}%
\special{pa 3280 3510}%
\special{fp}%
\special{pa 3280 2920}%
\special{pa 2890 3510}%
\special{fp}%
\special{pa 3680 2910}%
\special{pa 4080 3520}%
\special{fp}%
\special{pa 4090 2920}%
\special{pa 3680 3520}%
\special{fp}%
\special{pa 4480 2930}%
\special{pa 4870 3520}%
\special{fp}%
\special{pa 4880 2930}%
\special{pa 4490 3520}%
\special{fp}%
%
\special{pn 13}%
\special{pa 3500 840}%
\special{pa 4187 768}%
\special{fp}%
\special{pa 3726 1150}%
\special{pa 3959 473}%
\special{fp}%
\special{pa 3963 1475}%
\special{pa 4648 1418}%
\special{fp}%
\special{pa 4189 1784}%
\special{pa 4422 1108}%
\special{fp}%
\special{pa 4412 2108}%
\special{pa 5118 2046}%
\special{fp}%
\special{pa 4657 2426}%
\special{pa 4887 1728}%
\special{fp}%
\special{pa 4891 2731}%
\special{pa 5575 2672}%
\special{fp}%
\special{pa 5122 3047}%
\special{pa 5355 2371}%
\special{fp}%
%
\special{pn 13}%
\special{pa 2970 520}%
\special{pa 3242 1155}%
\special{fp}%
\special{pa 2741 827}%
\special{pa 3457 849}%
\special{fp}%
\special{pa 2501 1150}%
\special{pa 2758 1787}%
\special{fp}%
\special{pa 2272 1457}%
\special{pa 2987 1479}%
\special{fp}%
\special{pa 2029 1766}%
\special{pa 2297 2422}%
\special{fp}%
\special{pa 1798 2094}%
\special{pa 2533 2108}%
\special{fp}%
\special{pa 1576 2407}%
\special{pa 1835 3044}%
\special{fp}%
\special{pa 1342 2722}%
\special{pa 2057 2745}%
\special{fp}%
\put(56.2000,-27.5000){\makebox(0,0)[lb]{$E_m$}}%
\put(54.0000,-23.9000){\makebox(0,0)[lb]{$G_n$}}%
\put(51.6000,-21.1000){\makebox(0,0)[lb]{$E_n$}}%
\put(49.2000,-17.2000){\makebox(0,0)[lb]{$G_m$}}%
\put(47.0000,-14.8000){\makebox(0,0)[lb]{$F_{ik}$}}%
\put(44.4000,-11.4000){\makebox(0,0)[lb]{$F_{jl}$}}%
\put(42.2000,-8.2000){\makebox(0,0)[lb]{$F_{il}$}}%
\put(39.9000,-5.2000){\makebox(0,0)[lb]{$F_{jk}$}}%
\put(28.2000,-5.0000){\makebox(0,0)[lb]{$E_k$}}%
\put(25.4000,-9.1000){\makebox(0,0)[lb]{$G_l$}}%
\put(23.8000,-11.3000){\makebox(0,0)[lb]{$E_l$}}%
\put(20.6000,-15.3000){\makebox(0,0)[lb]{$G_k$}}%
\put(18.5000,-17.5000){\makebox(0,0)[lb]{$F_{im}$}}%
\put(15.4000,-21.5000){\makebox(0,0)[lb]{$F_{jn}$}}%
\put(14.6000,-23.8000){\makebox(0,0)[lb]{$F_{jm}$}}%
\put(11.4000,-28.1000){\makebox(0,0)[lb]{$F_{in}$}}%
\put(19.9000,-37.1000){\makebox(0,0)[lb]{$E_i$}}%
\put(23.9000,-37.3000){\makebox(0,0)[lb]{$G_j$}}%
\put(28.2000,-37.3000){\makebox(0,0)[lb]{$E_j$}}%
\put(32.1000,-37.2000){\makebox(0,0)[lb]{$G_i$}}%
\put(36.2000,-37.1000){\makebox(0,0)[lb]{$F_{km}$}}%
\put(40.3000,-37.1000){\makebox(0,0)[lb]{$F_{ln}$}}%
\put(44.4000,-37.1000){\makebox(0,0)[lb]{$F_{kn}$}}%
\put(48.4000,-37.2000){\makebox(0,0)[lb]{$F_{lm}$}}%
\end{picture}%
\end{center}
\caption{The $27$ lines on $\ol{\Sol}(\th)$ 
viewed from the tritangent plane at infinity} 
\label{fig:cubic3} 
\end{figure}

 
Moreover there are exactly $45$ tritangent planes that 
cut out a triplet of lines on $S$. 
In our case $S = \ol{\Sol}(\th)$, the plane at infinity 
$\{\, X \in \P^3 \,:\, X_0 = 0 \,\}$ is an instance of 
tritangent plane, which cuts out the lines in 
(\ref{eqn:lines}). 
Figure \ref{fig:cubic3} offers an arrangement of the $27$ lines 
viewed from the tritangent plane at infinity, 
where $\{i,j,k\} = \{1,2,3\}$ and $\{l,m,n\}=\{4,5,6\}$, and 
\begin{equation} \label{eqn:lineinf}
L_i = F_{ij},  \qquad L_j = F_{kl}, \qquad L_k = F_{mn} 
\end{equation}
are allocated for the lines at infinity. 
Each line at infinity is intersected by exactly eight lines 
and this fact enables us to divide the 27 lines into three 
groups of nine lines labeled by lines at infinity. 
{\sl Caution:} only the intersection relations among lines of 
the same group are indicated in Figure \ref{fig:cubic3}, with 
no other intersection relations being depicted. 
\par 
If $E_0$ is the strict transform of a plane in $\P^2$ not 
passing through $P_1, \dots, P_6$ 
relative to the $6$-point blow-up $S \to \P^2$, then the 
second cohomology group of $S = \ol{\Sol}(\th)$ is 
expressed as 
\begin{equation} \label{eqn:basis}
H^2(\ol{\Sol}(\th),\Z) = 
\Z E_0 \oplus \Z E_i \oplus \Z E_j \oplus \Z E_k \oplus 
\Z E_l \oplus \Z E_m \oplus \Z E_n, 
\end{equation}
where a divisor is identified with the cohomology class 
it represents. 
It is a Lorentz lattice of rank $7$ with intersection numbers 
\begin{equation} \label{eqn:in1}
(E_a, E_b) = \left\{ 
\begin{array}{ll}
{\-}1 \quad & (a = b = 0), \\[1mm] 
 -1 \quad & (a = b \neq 0), \\[1mm]
{\-}0 \quad & (\text{otherwise}). 
\end{array}\right. 
\end{equation}
In terms of the basis in (\ref{eqn:basis}) 
the lines $F_{ab}$ and $G_a$ are represented as 
\begin{equation} \label{eqn:FG}
F_{ab} = E_0 - E_a - E_b, \qquad 
G_a = 2 E_0 - (E_1 + \cdots + \widehat{E}_a + \cdots + E_6). 
\end{equation}
\par
We shall describe the $27$ lines on our cubic surface 
$\ol{\Sol}(\th)$ under the condition that $\ol{\Sol}(\th)$ 
is smooth, namely, $\th = \rh(\k)$ with $\k \in \K-\Wall$. 
To this end we introduce new parameters 
$b = (b_1,b_2,b_3,b_4) \in B := (\C_b^{\times})^4$ in such a 
manner that $b$ is expressed as 
\[ 
b_i = \left\{
\begin{array}{rl}
\exp(\sqrt{-1}\pi\k_i) \qquad  & (i = 1,2,3), \\[2mm]
-\exp(\sqrt{-1}\pi\k_4) \qquad & (i = 4), 
\end{array}
\right.
\]
as a function of $\k \in \K$. 
Then the Riemann scheme in Table \ref{tab:riemann} implies 
that $b_i$ is an eigenvalue of the monodromy matrix $\rho(\ga_i)$ 
around the point $t_i$ and formula (\ref{eqn:a}) implies that 
$a_i = b_i + b_i^{-1}$. 
Here parameters $b \in B$ should not be confused with the 
index $b$ above. 
In terms of the parameters $b \in B$, the discriminant $\vD(\th)$ 
of the cubic surfaces $\Sol(\th)$ factors as 
\begin{equation} \label{eqn:vD}
\vD(\th) = \prod_{l=1}^4(b_l-b_l^{-1})^2 
\prod_{\ve \in \{\pm1\}^4}(b^{\ve}-1), 
\end{equation}
where 
$b^{\ve} = b_1^{\ve_1}b_2^{\ve_2}b_3^{\ve_3}b_4^{\ve_4}$ 
for each quadruple sign 
$\ve = (\ve_1,\ve_2,\ve_3,\ve_4) \in \{\pm1\}^4$. 
Formula (\ref{eqn:vD}) clearly shows for which parameters 
$b \in B$ the cubic surface $\ol{\Sol}(\th)$ is smooth 
or singular. 
\par
Let $L_i(b_i,b_4;b_j,b_k)$ denote the line in $\P^3$ defined by 
the system of linear equations 
\begin{equation} \label{eqn:line2} 
X_i = (b_i b_4 + b_i^{-1}b_4^{-1}) X_0, \qquad 
X_j + (b_ib_4) X_k =  
\{ b_i (b_k+b_k^{-1}) + b_4 (b_j+b_j^{-1}) \} X_0. 
\end{equation}
Assume that $\ol{\Sol}(\th)$ is smooth, namely, $\vD(\th) \neq 0$. 
Then, as is mentioned earlier, for each $i \in \{1,2,3\}$ there 
are exactly eight lines on $\ol{\Sol}(\th)$ intersecting $L_i$, 
but not intersecting the remaining two lines at infinity, 
$L_j$ and $L_k$. 
They are just $\{E_i, G_j\}$, $\{E_j, G_i\}$, $\{F_{km}, F_{ln}\}$, 
$\{F_{kn}, F_{lm}\}$ as in Figure \ref{fig:cubic3}, where 
two lines from the same pair intersect, but ones from 
different pairs are disjoint. 
In terms of parameters $b \in B$ those eight lines are 
given as in Table \ref{tab:line}. 
\begin{table}[t] 
\begin{center} 
\begin{tabular}{|c|c|c|}
\hline
\vspace{-0.3cm} & & \\
$1$ & $L_i(b_i,b_4;b_j,b_k)$ & $L_i(1/b_i,1/b_4;b_j,b_k)$ \\[2mm]
\hline
\hline
\vspace{-0.3cm} & & \\
$2$ & $L_i(b_j,b_k;b_i,b_4)$ & $L_i(1/b_j,1/b_k;b_i,b_4)$ \\[2mm]
\hline
\hline
\vspace{-0.3cm} & & \\
$3$ & $L_i(1/b_i,b_4;b_j,b_k)$ & $L_i(b_i,1/b_4;b_j,b_k)$ \\[2mm]
\hline
\hline
\vspace{-0.3cm} & & \\
$4$ & $L_i(1/b_j,b_k;b_i,b_4)$ & $L_i(b_j,1/b_k;b_i,b_4)$ \\[2mm]
\hline
\end{tabular}
\end{center} 
\caption{Eight lines intersecting the line $L_i$ at infinity, 
divided into four pairs} 
\label{tab:line}
\end{table}
\section{Involutions on Cubic Surface} \label{sec:involution}
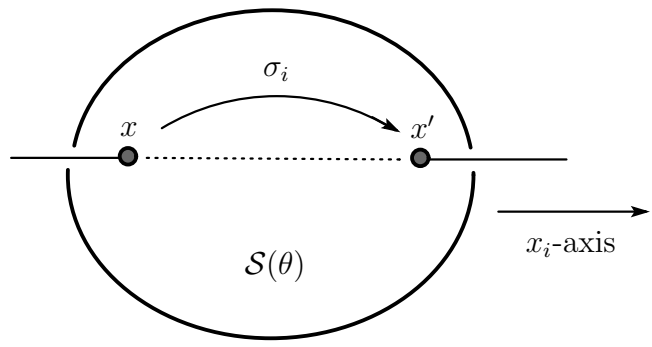
\begin{figure}[t]
\begin{center}
\unitlength 0.1in
\begin{picture}(33.24,17.19)(4.90,-19.59)
%
\special{pn 13}%
\special{pa 490 1014}%
\special{pa 1111 1014}%
\special{fp}%
%
\special{pn 13}%
\special{pa 2641 1023}%
\special{pa 3397 1023}%
\special{fp}%
%
\special{pn 20}%
\special{ar 1858 1095 1053 855  3.3124258 6.1241532}%
%
\special{pn 20}%
\special{ar 1849 1104 1062 855  6.2694533 6.2831853}%
\special{ar 1849 1104 1062 855  0.0000000 3.1735968}%
\put(17.1400,-16.6200){\makebox(0,0)[lb]{$\Sol(\th)$}}%
%
\special{pn 13}%
\special{pa 1201 1014}%
\special{pa 2524 1023}%
\special{dt 0.045}%
\special{pa 2524 1023}%
\special{pa 2523 1023}%
\special{dt 0.045}%
\put(10.6000,-9.0000){\makebox(0,0)[lb]{$x$}}%
\put(25.8000,-9.1000){\makebox(0,0)[lb]{$x'$}}%
%
\special{pn 20}%
\special{sh 0.600}%
\special{ar 1102 1005 44 44  0.0000000 6.2831853}%
%
\special{pn 20}%
\special{sh 0.600}%
\special{ar 2632 1014 44 44  0.0000000 6.2831853}%
\put(18.1000,-6.1000){\makebox(0,0)[lb]{$\si_i$}}%
%
\special{pn 13}%
\special{pa 3040 1295}%
\special{pa 3814 1295}%
\special{fp}%
\special{sh 1}%
\special{pa 3814 1295}%
\special{pa 3747 1275}%
\special{pa 3761 1295}%
\special{pa 3747 1315}%
\special{pa 3814 1295}%
\special{fp}%
\put(31.8400,-15.2900){\makebox(0,0)[lb]{$x_i$-axis}}%
%
\special{pn 13}%
\special{ar 1880 1850 1161 1161  4.1675251 5.2764727}%
%
\special{pn 13}%
\special{pa 2380 800}%
\special{pa 2510 870}%
\special{fp}%
\special{sh 1}%
\special{pa 2510 870}%
\special{pa 2461 821}%
\special{pa 2463 845}%
\special{pa 2442 856}%
\special{pa 2510 870}%
\special{fp}%
\end{picture}%
\end{center}
\caption{Involutions of a $(2,2,2)$-surface} 
\label{fig:involution}
\end{figure}

The affine cubic surface $\Sol(\th)$ is a $(2,2,2)$-surface, 
that is, its defining equation $f(x,\th) = 0$ is a quadratic 
equation in each variable $x_i$, $i = 1,2,3$. 
Therefore the line through a point $x \in \Sol(\th)$ parallel 
to the $x_i$-axis passes through a unique second point 
$x' \in \Sol(\th)$ (see Figure \ref{fig:involution}). 
This defines an involution $\si_i : \Sol(\th) \to \Sol(\th)$, 
$x \mapsto x'$, which is explicitly given by 
\begin{equation} \label{eqn:si} 
\si_i : \qquad 
(x_i',x_j', x_k') = (\th_i-x_i-x_jx_k, x_j, x_k), 
\qquad (i = 1,2,3).  
\end{equation}
\par
\begin{figure}[t] 
\begin{center}
\unitlength 0.1in
\begin{picture}(24.70,21.70)(6.00,-24.70)
\put(17.3000,-24.1000){\makebox(0,0)[lb]{$L_i$}}%
\put(9.6000,-17.4000){\makebox(0,0)[lb]{$L_j$}}%
\put(25.1000,-17.4000){\makebox(0,0)[lb]{$L_k$}}%
\put(28.6000,-23.6000){\makebox(0,0)[lb]{$p_j$}}%
\put(17.5000,-4.7000){\makebox(0,0)[lb]{$p_i$}}%
\put(6.2000,-23.5000){\makebox(0,0)[lb]{$p_k$}}%
\put(16.2000,-17.5000){\makebox(0,0)[lb]{$\si_i \car$}}%
\put(22.6000,-13.1000){\makebox(0,0)[lb]{$q_k$}}%
\put(12.2000,-13.2000){\makebox(0,0)[lb]{$q_j$}}%
%
\special{pn 20}%
\special{sh 0.600}%
\special{ar 1820 630 40 40  0.0000000 6.2831853}%
%
\special{pn 20}%
\special{sh 0.600}%
\special{ar 2720 2160 40 40  0.0000000 6.2831853}%
%
\special{pn 20}%
\special{sh 0.600}%
\special{ar 920 2160 40 40  0.0000000 6.2831853}%
%
\special{pn 20}%
\special{pa 1640 330}%
\special{pa 2900 2460}%
\special{fp}%
%
\special{pn 20}%
\special{pa 1980 330}%
\special{pa 740 2470}%
\special{fp}%
%
\special{pn 20}%
\special{pa 600 2170}%
\special{pa 3070 2170}%
\special{fp}%
%
\special{pn 20}%
\special{sh 0.600}%
\special{ar 1370 1390 40 40  0.0000000 6.2831853}%
%
\special{pn 20}%
\special{sh 0.600}%
\special{ar 2260 1390 40 40  0.0000000 6.2831853}%
%
\special{pn 13}%
\special{ar 1830 1660 1015 1015  4.9362935 6.2831853}%
\special{ar 1830 1660 1015 1015  0.0000000 0.2970642}%
%
\special{pn 13}%
\special{ar 1820 1650 1016 1016  2.8289604 4.4722180}%
%
\special{pn 13}%
\special{pa 2820 1870}%
\special{pa 2810 1960}%
\special{fp}%
\special{sh 1}%
\special{pa 2810 1960}%
\special{pa 2837 1896}%
\special{pa 2816 1907}%
\special{pa 2797 1892}%
\special{pa 2810 1960}%
\special{fp}%
%
\special{pn 13}%
\special{pa 820 1840}%
\special{pa 840 1960}%
\special{fp}%
\special{sh 1}%
\special{pa 840 1960}%
\special{pa 849 1891}%
\special{pa 831 1907}%
\special{pa 809 1898}%
\special{pa 840 1960}%
\special{fp}%
%
\special{pn 13}%
\special{pa 2180 700}%
\special{pa 2070 670}%
\special{fp}%
\special{sh 1}%
\special{pa 2070 670}%
\special{pa 2129 707}%
\special{pa 2121 684}%
\special{pa 2140 668}%
\special{pa 2070 670}%
\special{fp}%
%
\special{pn 13}%
\special{pa 1450 700}%
\special{pa 1570 660}%
\special{fp}%
\special{sh 1}%
\special{pa 1570 660}%
\special{pa 1500 662}%
\special{pa 1519 677}%
\special{pa 1513 700}%
\special{pa 1570 660}%
\special{fp}%
%
\special{pn 13}%
\special{pa 1430 2080}%
\special{pa 1790 780}%
\special{fp}%
\special{sh 1}%
\special{pa 1790 780}%
\special{pa 1753 839}%
\special{pa 1776 831}%
\special{pa 1791 850}%
\special{pa 1790 780}%
\special{fp}%
%
\special{pn 13}%
\special{pa 2170 2080}%
\special{pa 1850 780}%
\special{fp}%
\special{sh 1}%
\special{pa 1850 780}%
\special{pa 1847 850}%
\special{pa 1863 832}%
\special{pa 1885 840}%
\special{pa 1850 780}%
\special{fp}%
\end{picture}%
\end{center}
\caption{The birational map $\si_i$ restricted to $L$}
\label{fig:cubic2} 
\end{figure}
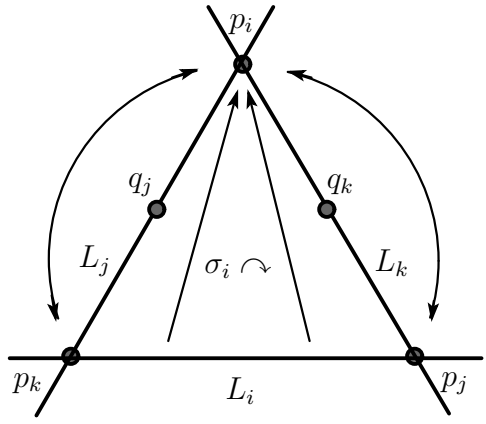

 
The automorphism $\si_i$ of the affine surface 
$\Sol(\th)$ extends to a birational map of the projective surface 
$\ol{\Sol}(\th)$, which will also be denoted by $\si_i$. 
In terms of the homogeneous coordinates $X$ of $\P^3$, the 
birational map $\si_i : X \mapsto X'$ is expressed as 
\[
[X_0':X_i':X_j':X_k'] = 
[X_0^2 : \th_i X_0^2 - X_0 X_i - X_j X_k : X_0 X_j : X_0 X_k] 
\]
We shall investigate the behavior of the birational map $\si_i$ 
in a neighborhood of the tritangent lines $L$ at infinity. 
To this end let us introduce the following three points 
\[
q_1 = [0:0:1:1], \qquad 
q_2 = [0:1:0:1], \qquad 
q_3 = [0:1:1:0],  
\]
where $q_i$ may be thought of as the ``mid-point" of $p_j$ and $p_k$ 
on the line $L_i$. 
\begin{lemma} \label{lem:L} 
The birational map $\si_i$ has the following properties 
$($see Figure $\ref{fig:cubic2}$$)$. 
\begin{enumerate}
\item $\si_i$ blows down the line $L_i$ to the point $p_i$,  
\item $\si_i$ restricts to the automorphism of $L_j$ that fixes 
$q_j$ and exchanges $p_i$ and $p_k$, 
\item $\si_i$ restricts to the automorphism 
of $L_k$ that fixes $q_k$ and exchanges $p_i$ and $p_j$, 
\item $p_i$ is the unique indeterminacy point of $\si_i$, 
\end{enumerate}
\end{lemma} 
{\it Proof}. 
In order to investigate $\si_i$, we make use of inhomogeneous 
coordinates of $\P^3$ in (\ref{eqn:coord}) and local 
coordinates and local equations of $\ol{\Sol}(\th)$ in 
Table \ref{tab:local}, with target coordinates 
being dashed. 
\par 
In terms of inhomogeneous coordinates $v$ and $u'$ of $\P^3$, 
the map $\si_i : v \mapsto u'$ is expressed as 
\begin{equation} \label{eqn:vu'}
u_0' = \frac{v_0^2}{\th_i v_0^2 - v_0 v_i - v_k}, \qquad  
u_j' = \frac{v_0}{\th_i v_0^2 - v_0 v_i - v_k},   \qquad 
u_k' = \frac{v_0 v_k}{\th_i v_0^2 - v_0 v_i - v_k}. 
\end{equation}
In a neighborhood of $L_i-\{p_j,p_k\}$ in $\ol{\Sol}(\th)$, 
using $v_i = O(v_0)$, we observe that 
\[
\th_i v_0^2-v_0 v_i-v_k = -v_k \{1+O(v_0^2)\}, 
\]
which is substituted into (\ref{eqn:vu'}) to yield 
\[
u_j' = -\frac{v_0}{v_k \{1+O(v_0^2)\}} = 
-\frac{v_0}{v_k}\{1 + O(v_0^2)\}, 
\qquad 
u_k' = -\frac{v_0 v_k}{v_k \{1+O(v_0^2)\}} = -v_0\{1 + O(v_0^2)\}.
\]
In particular putting $v_0 = 0$ leads to $u_j' = u_k' = 0$. 
This means that $\si_i$ maps a neighborhood of $L_i - \{p_j,p_k \}$ 
to a neighborhood of $p_i$, collapsing $L_i - \{p_j,p_k \}$ 
to the single point $p_i$. 
\par 
In a similar manner, in a neighborhood of $p_j$ in $\ol{\Sol}(\th)$ 
we observe that 
\[
v_0 = -(v_i v_k) \{1+O_2(v_i,v_k)\}, \qquad 
\th_i v_0^2-v_0 v_i-v_k = -v_k \{1+O_2(v_i, v_k)\}, 
\] 
which are substituted into (\ref{eqn:vu'}) to yield 
\[
u_j' = v_i \{1+O_2(v_i,v_k)\}, 
\qquad 
u_k' = (v_i v_k) \{1+O_2(v_i,v_k)\}. 
\]
In particular putting $v_i = 0$ leads to $u_j' = u_k' = 0$. 
This means that $\si_i$ maps a neighborhood of $p_j$ to 
a neighborhood of $p_i$, collapsing a neighborhood in $L_i$ 
of $p_j$ to the single point $p_i$. 
Using $w$ and $u'$ in place of $v$ and  $u'$, we can 
make a similar argument in a neighborhood of $p_k$.  
Therefore $\si_i$ blows down $L_i$ to the point $p_i$, which 
proves assertion (1). 
Moreover it is clear from the argument that there is no 
indeterminacy point on the line $L_i$. 
\par 
In terms of inhomogeneous coordinates $u$ and $u'$ of 
$\P^3$ the map $\si_i : u \mapsto u'$ is expressed as  
\begin{equation} \label{eqn:siu}
u_0' = \frac{u_0^2}{\th_i u_0^2 - u_0 -u_j u_k},   \quad  
u_j' = \frac{u_0 u_j}{\th_i u_0^2 - u_0 -u_j u_k}, \quad 
u_k' = \frac{u_0 u_k}{\th_i u_0^2 - u_0 -u_j u_k}.
\end{equation}
In a neighborhood of $L_j - \{ p_i,p_k \}$ in $\ol{\Sol}(\th)$, 
using $u_j = -(u_k + 1/u_k) u_0 + O(u_0^2)$, we have 
\[
\th_i u_0^2 - u_0 -u_j u_k = u_0 \{ u_k^2 + O(u_0) \},
\] 
which is substituted into (\ref{eqn:siu}) to yield 
\[
u_0' = \frac{u_0}{u_k^2 + O(u_0)} = \frac{u_0}{u_k^2} + O(u_0^2), 
\qquad 
u_k' = \frac{u_k}{u_k^2 + O(u_0)} = \frac{1}{u_k} + O(u_0).
\]
In particular putting $u_0 = 0$ leads to $u_0' = 0$ and 
$u_k' = 1/u_k$. 
This means that $\si_i$ restricts to an automorphism of 
a neighborhood of $L_j - \{p_i,p_k\}$ in $\ol{\Sol}(\th)$ 
that induces a unique automorphism of $L_i$ fixing 
$q_j$ and exchanging $p_i$ and $p_k$. 
This proves assertion (2) and also shows that 
there is no indeterminacy point on $L_j-\{p_i,p_k\}$. 
Assertion (3) and the nonexistence of indeterminacy point 
on $L_k-\{p_i,p_j\}$ are established just in the same manner. 
\par 
From the above argument we have already known that 
there is no indeterminacy point other than $p_i$. 
Then the point $p_i$ is actually an indeterminacy point, 
because $\si_i$ is an involution blowing down 
$L_i$ to $p_i$ and hence blows up $p_i$ to $L_i$ reciprocally. 
This proves assertion (4).  
\hfill $\Box$ \par\medskip 
Later we will need some information about how the involution $\si_i$ 
transforms a line to another curve, which is stated in the following 
lemma. 
\begin{lemma} \label{lem:si} 
The involution $\si_i$ satisfies the following properties: 
\begin{enumerate} 
\item $\si_i(E_i)$ intersects $E_i$ at two points counted 
with multiplicity, 
\item $\si_i(E_i)$ intersects $E_j$ at one point counted with 
multiplicity, 
\item $\si_i$ exchanges the lines $E_k$ and $G_l$; $E_l$ and $G_k$; 
$E_m$ and $G_n$; $E_n$ and $G_m$, respectively. 
\end{enumerate}
\end{lemma} 
{\it Proof}. 
By Table \ref{tab:line} we may put $E_i = L_i(b_i,b_4;b_j,b_k)$ 
and $E_j = L_i(b_j,b_k;b_i,b_4)$. 
Assertion (1) of Lemma \ref{lem:L} implies that $\si_i(E_i)$ 
does not intersect $E_i$ nor $E_j$ at any point at infinity. 
So we can work with inhomogeneous coordinates 
$x = (x_1,x_2,x_3)$. 
In view of (\ref{eqn:line2}) the line $E_i$ is given by 
\begin{equation} \label{eqn:Ei} 
x_i = b_ib_4 + (b_ib_4)^{-1}, \qquad 
x_j + (b_ib_4) x_k =  a_k b_i +  a_j b_4. 
\end{equation}
In a similar manner, by exchanging $(b_i,b_4)$ and $(b_j,b_k)$ in 
(\ref{eqn:line2}), the line $E_j$ is given by 
\begin{equation} \label{eqn:Ej} 
x_i = b_jb_k + (b_jb_k)^{-1}, \qquad 
x_j + (b_jb_k) x_k =  a_4 b_j +  a_i b_k. 
\end{equation}
Moreover, by applying formula (\ref{eqn:si}) to (\ref{eqn:Ei}), 
the curve $\si_i(E_i)$ is expressed as 
\begin{equation} \label{eqn:si(Ei)}
\th_i-x_i-x_jx_k = b_ib_4+(b_ib_4)^{-1}, \qquad 
x_j + (b_ib_4) x_k =  a_k b_i +  a_j b_4. 
\end{equation}
Note that the second equations of 
(\ref{eqn:Ei}) and (\ref{eqn:si(Ei)}) are the same. 
\par 
In order to find out the intersection of $\si_i(E_i)$ with 
$E_i$, let us couple (\ref{eqn:Ei}) and (\ref{eqn:si(Ei)}). 
Eliminating $x_i$ and $x_j$ we obtain 
a quadratic equation for $x_k$, 
\[
(b_ib_4) x_k^2 -(a_k b_i + a_j b_4)x_k + 
\th_i-2\{b_ib_4 + (b_ib_4)^{-1}\} = 0. 
\]
For each root of this equation we have an intersection point 
of $\si_i(E_i)$ with $E_i$; for a double root we have an 
intersection point of multiplicity two. 
This proves assertion (1). 
\par 
Next, in order to find out the intersection of $\si_i(E_i)$ 
with $E_j$, let us couple (\ref{eqn:Ej}) and (\ref{eqn:si(Ei)}). 
From the first equation of (\ref{eqn:Ej}) the 
$x_i$-coordinate is already fixed. 
The second equations of (\ref{eqn:Ej}) and (\ref{eqn:si(Ei)}) 
are coupled to yield a linear system for $x_j$ and $x_k$, 
whose determinant 
\[
b_jb_k-b_ib_4= b_ib_4(b_i^{-1}b_jb_kb_4^{-1}-1)
\]
is nonzero from the assumption that $\ol{\Sol}(\th)$ is smooth, 
that is, the discriminant $\vD(\th)$ in (\ref{eqn:vD}) is nonzero. 
Then the linear system is uniquely solved to determine 
$x_j$ and $x_k$. 
Now we can check that the first equation of (\ref{eqn:si(Ei)}) 
is redundant, that is, automatically satisfied. 
Therefore $\si_i(E_i)$ and $E_j$ has a simple intersection, 
which implies assertion (2). 
\par 
Finally we see that $\si_i$ exchanges $E_k$ and $G_l$. 
We may put $E_k = L_j(b_j,b_4;b_k,b_i)$ and 
$G_l = L_j(1/b_j,1/b_4;b_k,b_i)$. 
By formula (\ref{eqn:line2}) (with indices suitably permuted), 
these lines are given by 
\begin{align} 
x_j &= b_jb_4 + (b_jb_4)^{-1}, & 
x_k + (b_jb_4) x_i &= a_ib_j + a_kb_4, \label{eqn:Ek} 
\\[2mm] 
x_j &= b_jb_4 + (b_jb_4)^{-1}, &
x_k + (b_jb_4)^{-1} x_i &= a_i b_j^{-1} + a_k b_4^{-1}. 
\label{eqn:Gl} 
\end{align} 
Using formula (\ref{eqn:si}) we can check that equations 
(\ref{eqn:Ek}) and (\ref{eqn:Gl}) are transformed to each other 
by $\si_i$. 
This together with similar argument for the other lines 
establishes assertion (3). \hfill $\Box$ 
\section{Dynamical System on Cubic Surface} \label{sec:dynamics} 
Let $G = \la \si_1, \si_2, \si_3 \ra$ be the group of birational 
transformations on $\ol{\Sol}(\th)$ generated by the involutions 
$\si_1$, $\si_2$, $\si_3$. 
We are interested in the dynamics of the $G$-action on 
$\ol{\Sol}(\th)$. 
Usually the dynamics of a group action is more involved than 
that of a single transformation; more techniques 
and tools have been developed for the latter rather than for 
the former. 
So it may be better to pick up a single transformation from the 
group $G$ and study its dynamics. 
For such a transformation we take a composition of the three 
basic involutions, 
\begin{equation} \label{eqn:coxeter}
c = \si_i \circ \si_j \circ \si_k \, : \,
\ol{\Sol}(\th) \carl.
\end{equation}
If $G$ is regarded as a nonlinear reflection group with basic 
`reflections' $\si_1$, $\si_2$, $\si_3$, then $c$ may be 
thought of as a `Coxeter' transformation and it is expected 
that the dynamics of the transformation $c$ plays a dominant 
role in understanding the dynamics of the whole $G$-action. 
\par
The relevance of the transformation (\ref{eqn:coxeter}) to 
our main problem is stated as follows. 
\begin{lemma} \label{lem:coxeter} 
Via the Riemann-Hilbert correspondence $(\ref{eqn:RHtk2})$ 
the Pochhammer-Poincar\'e map $\wp_* : M_x(\k) \carl$ is 
strictly conjugated to the square $c^2 : \Sol(\th) \carl$ 
of the Coxeter transformation $(\ref{eqn:coxeter})$, 
restricted to the affine part $\Sol(\th)$ of the cubic 
surface $\ol{\Sol}(\th)$. 
\end{lemma}
{\it Proof}. 
Since the transformation $g_i$ in (\ref{eqn:gi}) is a strict 
conjugacy of the half-Poincar\'e map $\b_{i*}$ in 
(\ref{eqn:PSi}), a glance at (\ref{eqn:commutator1}) and 
(\ref{eqn:PSi}) shows that the commutator 
$[g_i^2, g_j^{-2}] = g_i^2 g_j^{-2} g_i^{-2} g_j^2$ is a 
strict conjugacy of the Pochhammer-Poincar\'e map $\wp_*$. 
On the other hand, using (\ref{eqn:gi}) and 
(\ref{eqn:si}), we can directly check that 
$g_i^2 g_j^{-2} g_i^{-2} g_j^2 =(\si_i\si_j\si_k)^2=c^2$. 
Hence $c^2$ is a strict conjugacy of $\wp_*$. 
\hfill $\Box$ \par\medskip
A general theory of dynamical systems for bimeromorphic maps of 
surfaces is developed in \cite{DF}. 
We shall apply it to our map (\ref{eqn:coxeter}) upon 
reviewing some rudiments of the article \cite{DF}. 
Let $S$ be a compact complex surface, $f : S \carl$ 
a bimeromorphic map. 
Then $f$ is represented by a compact complex surface $\vG$ 
and proper modifications $\pi_1 : \vG \to S$ and 
$\pi_2 : \vG \to S$ such that $f = \pi_2 \ci \pi_1^{-1}$ on a 
dense open subset. 
For $i = 1, 2$, let $\E(\pi_i) := \{\, x \in \vG \,:\, 
\mathrm{\#}    \, \pi_i^{-1}(\pi_i (x)) = \infty \,\}$ 
be the exceptional set for the projection $\pi_i$. 
The images $I(f) := \pi_1 (\E(\pi_1))$ and 
$\E(f) := \pi_1 (\E(\pi_2))$ are called the indeterminacy 
set and the exceptional set of $f$ respectively. 
In our case where $S = \ol{\Sol}(\th)$ and $f = \si_i$, 
Lemma \ref{lem:L} implies that these sets are described 
as follows. 
\begin{lemma} \label{lem:IE}
$I(\si_i) = \{p_i\}$ and $\E(\si_i) = L_i$ for $i = 1, 2, 3$.
\end{lemma} 
\par 
If $S$ is K\"ahler, two natural actions of $f$, 
pull-back and push-forward, on the Dolbeault 
cohomology group $H^{1,1}(S)$ are defined in the following 
manner: 
A smooth $(1,1)$-form $\omega$ on $S$ can be pulled back as 
a smooth $(1,1)$-form $\pi_2^* \omega$ on $\vG$ and then 
pushed forward as a $(1,1)$-current $\pi_{1*}\pi_2^* \omega$ 
on $S$. 
Hence we define the pull-back 
$f^* \omega := \pi_{1*}\pi_2^* \omega$ and also the 
push-forward 
$f_* \omega = (f^{-1})^* \omega := \pi_{2*}\pi_1^* \omega$. 
The operators $f^*$ and $f_*$ commute with the exterior 
differential $d$ and the complex structure of $S$ and 
so descend to linear actions on $H^{1,1}(S)$. 
For general bimeromorphic maps $f$ and $g$, the 
composition rule $(f\ci g)^* = g^* \ci f^*$ is not 
necessarily true. 
But a useful criterion under which this rule becomes 
true is given in \cite{DF}. 
\begin{lemma} \label{lem:compo}
If $f(\E(f)) \cap I(g) = \emptyset$, then 
$(f \ci g)^* = g^* \ci f^* : H^{1,1}(S) \carl$.
\end{lemma}
\par 
We apply this lemma to our Coxeter transformation 
$c = \si_i \ci \si_j \ci \si_k$. 
\begin{lemma} \label{lem:compo2}
We have $c^{*} = \si_k^{*} \ci \si_j^{*} \ci \si_i^{*} 
: H^{1,1}(\ol{\Sol}(\th)) \carl$.
\end{lemma}
{\it Proof}.  
First we apply Lemma \ref{lem:compo} to 
$f = \si_i$ and $g = \si_j \ci \si_k$. 
By Lemmas \ref{lem:L} and \ref{lem:IE} we have 
$\E(\si_i) = L_i$ and $I(\si_j \ci \si_k) = \{p_k\}$ and 
so $\si_i(\E(\si_i)) \cap I(\si_j \ci \si_k) = 
\{p_i \} \cap \{ p_k \} = \emptyset$, which means that 
the condition of Lemma \ref{lem:compo} is satisfied. 
Then the lemma yields $(\si_i \ci \si_j \ci \si_k)^* 
= (\si_j \ci \si_k)^* \ci \si_i^*$. 
Next we apply Lemma \ref{lem:compo} to $f = \si_j$ 
and $g = \si_k$. 
Again by Lemmas \ref{lem:L} and \ref{lem:IE} we have 
$\E(\si_j) = L_j$ and $I(\si_k) = \{p_k \}$ and so 
$\si_j(\E(\si_j)) \cap I(\si_k) = \{p_j\} \cap \{p_k\} 
= \emptyset$, which means that the condition of Lemma 
\ref{lem:compo} is satisfied. 
Then the lemma yields $(\si_j \ci \si_k)^* = 
\si_k^* \ci \si_j^*$. 
Putting these two steps together, we obtain 
$c^* = (\si_i \ci \si_j \ci \si_k)^* = 
(\si_j \ci \si_k)^* \ci \si_i^* 
= \si_k^* \ci \si_j^* \ci \si_i^*$. 
\hfill $\Box$ \par\medskip
By Lemma \ref{lem:compo2} the calculation of the action 
$c^{*} : H^{1,1}(\ol{\Sol}(\th)) \carl$ is reduced to 
that of the actions $\si_i^*$, $\si_j^*$, 
$\si_k^* : H^{1,1}(\ol{\Sol}(\th)) \carl$, 
which is now set forth. 
Since the cubic surface $\ol{\Sol}(\th)$ is rational, 
we have $H^{1,1}(\ol{\Sol}(\th)) = H^2(\ol{\Sol}(\th),\C)$, 
where the latter group is described in (\ref{eqn:basis}). 
\begin{table}[t]
\[
\begin{array}{rl}
\si_i^* = \begin{pmatrix}
                {\-} 6 & {\-} 3 & {\-} 3 & {\-} 2 & {\-} 2 & 
                {\-} 2 & {\-} 2 \\ 
                -3 & -2 & -1 & -1 & -1 & -1 & -1 \\ 
                -3 & -1 & -2 & -1 & -1 & -1 & -1 \\ 
                -2 & -1 & -1 & -1 & {\-} 0 & -1 & -1 \\ 
                -2 & -1 & -1 & {\-} 0 & -1 & -1 & -1 \\ 
                -2 & -1 & -1 & -1 & -1 & -1 & {\-} 0 \\ 
                -2 & -1 & -1 & -1 & -1 & {\-} 0 & -1 
            \end{pmatrix} 
\quad & 
\si_j^* = \begin{pmatrix}
                {\-} 6 & {\-} 2 & {\-} 2 & {\-} 3 & {\-} 3 & 
                {\-} 2 & {\-} 2 \\
                -2 & -1 & {\-} 0 & -1 & -1 & -1 & -1 \\ 
                -2 & {\-} 0 & -1 & -1 & -1 & -1 & -1 \\ 
                -3 & -1 & -1 & -2 & -1 & -1 & -1 \\ 
                -3 & -1 & -1 & -1 & -2 & -1 & -1 \\ 
                -2 & -1 & -1 & -1 & -1 & -1 & {\-} 0 \\
                -2 & -1 & -1 & -1 & -1 & {\-} 0 & -1  
            \end{pmatrix} 
\\ &  \\
\si_k^* = \begin{pmatrix}
               {\-} 6 & {\-} 2 & {\-} 2 & {\-} 2 & {\-} 2 & 
               {\-} 3 & {\-} 3 \\ 
               -2 & -1 & {\-} 0 & -1 & -1 & -1 & -1 \\
               -2 & {\-} 0 & -1 & -1 & -1 & -1 & -1 \\
               -2 & -1 & -1 & -1 & {\-} 0 & -1 & -1 \\
               -2 & -1 & -1 & {\-} 0 & -1 & -1 & -1 \\
               -3 & -1 & -1 & -1 & -1 & -2 & -1 \\
               -3 & -1 & -1 & -1 & -1 & -1 & -2  
            \end{pmatrix}
\quad & 
c^* = \begin{pmatrix}
               12 & {\-} 6 & {\-} 6 & {\-} 4 & {\-} 4 & 
               {\-} 3 & {\-} 3 \\
               -3 & -2 & -1 & -1 & -1 & -1 & -1 \\ 
               -3 & -1 & -2 & -1 & -1 & -1 & -1 \\ 
               -4 & -2 & -2 & -2 & -1 & -1 & -1 \\ 
               -4 & -2 & -2 & -1 & -2 & -1 & -1 \\ 
               -6 & -3 & -3 & -2 & -2 & -2 & -1 \\ 
               -6 & -3 & -3 & -2 & -2 & -1 & -2  
        \end{pmatrix}
\end{array}
\]
\caption{Matrix representations of $\si_i^*$, $\si_j^*$, 
$\si_k^*$, $c^* : H^2(\ol{\Sol}(\th), \Z) \carl$}
\label{tab:matrix} 
\end{table}
\begin{proposition} \label{prop:in2} 
The linear maps $\si_i^*$, $\si_j^*$, $\si_k^*$,  
$c^* : H^2(\ol{\Sol}(\th), \Z) \carl$ admit matrix 
representations as in Table $\ref{tab:matrix}$ 
with respect to the basis in $(\ref{eqn:basis})$.  
The characteristic polynomial of $c^*$ is given by 
\begin{equation} \label{eqn:char}
\det(x I - c^*) = x(x+1)^4 (x^2 - 4 x -1), 
\end{equation} 
and hence its eigenvalues are $0$, $-1$ and 
$2 \pm \sqrt{5}$, where the eigenvalue $-1$ is quadruple 
while the remaining ones are all simple. 
The spectral radius $\rho(c^*)$ of $c^*$ is given by 
$2 + \sqrt{5}$. 
\end{proposition} 
{\it Proof}. 
First we shall find the matrix representation of $\si_i^*$. 
If $\xi_{ab}$ denotes the $(a,b)$-th entry of the matrix to 
be found, then (\ref{eqn:in1}) implies that 
\[
\si_i^* E_b = \sum_{a=0}^6 \xi_{ab} \, E_a 
= \sum_{a=0}^6 \delta_a (\si_i^* E_b, E_a) \, E_a, 
\]
where we put $\delta_a = 1$ for $a = 0$ and $\delta_a = -1$ for 
$a \neq 0$. 
Now we claim that 
\begin{equation} \label{eqn:xi}
\xi_{ab} = \delta_a (\si_i^*E_b, E_a), \qquad 
\xi_{ab} = \delta_a \delta_b \xi_{ba}. 
\end{equation}
The first formula in (\ref{eqn:xi}) is obvious and the second 
formula is derived as follows: 
\[
\xi_{ab} = \delta_a (\si_i^*E_b, E_a) = 
\delta_a (E_b, \si_{i*}E_a) = \delta_a (E_b, \si_i^*E_a) 
= (\delta_a \delta_b) \cdot \delta_b (\si_i^*E_a, E_b) = 
(\delta_a \delta_b) \xi_{ba}, 
\]
where in the third equality we have used the fact that 
$\si_i$ is an involution; $\si_{i*} = (\si_i^{-1})^* = 
\si_i^*$. 
By assertions (1) and (2) of Lemma \ref{lem:si} we have 
$(\si_i^*E_i,E_i) = 2$ and $(\si_i^*E_i, E_j) = 1$ and 
likewise $(\si_i^*E_j,E_j) = 2$ and 
$(\si_i^*E_j, E_i) = 1$. 
Then the first formula of (\ref{eqn:xi}) yields 
\begin{equation} \label{eqn:xi2} 
\xi_{ii} = \xi_{jj} = -2, \qquad 
\xi_{ij} = \xi_{ji} = -1. 
\end{equation}
The assertion (3) of Lemma \ref{lem:si} together 
with the second formula of (\ref{eqn:FG}) yields 
\begin{equation} \label{eqn:si3} 
\left\{
\begin{array}{rcl}
\si_i^*E_k &=&2 E_0-E_i-E_j-E_k \phantom{- E_l\,}-E_m-E_n,   \\[1mm]
\si_i^*E_l &=&2 E_0-E_i-E_j \phantom{- E_k\,\,}-E_l-E_m-E_n, \\[1mm]
\si_i^*E_m &=&2 E_0-E_i-E_j-E_k-E_l-E_m \phantom{- E_n\,},   \\[1mm]
\si_i^*E_n &=&2 E_0-E_i-E_j-E_k-E_l \phantom{- E_m\,}-E_n,   \\
\end{array}
\right.
\end{equation}
It follows from (\ref{eqn:xi2}) and (\ref{eqn:si3}) that the 
matrix representation for $\si_i^*$ takes the form 
\begin{equation} \label{eqn:si4} 
\si_i^* = 
\left(
\begin{array}{ccc|cccc}
    * & *  &  * & {\-}2 & {\-}2 & {\-}2 & {\-}2  \\
    * & -2 & -1 & -1    & -1    & -1    & -1     \\
    * & -1 & -2 & -1    & -1    & -1    & -1     \\[1mm]
\hline
    * & *  & *  & -1    & {\-}0 & -1    & -1     \\
    * & *  & *  & {\-}0 & -1    & -1    & -1     \\
    * & *  & *  & -1    & -1    & -1    & {\-} 0 \\
    * & *  & *  & -1    & -1    & {\-}0 & -1 
\end{array}
\right), 
\end{equation}
where the entries denoted by $*$ are yet to be determined. 
The $(2,1)$-block of (\ref{eqn:si4}) is easily determined 
by the second formula in (\ref{eqn:xi}). 
The final ingredient taken into account is the fact that 
$\si_i$ blows down $L_i = E_0-E_i-E_j$ to a point $p_i$ 
(see Lemma \ref{lem:L}), which leads to 
\[
\si_i^* E_0 -\si_i^* E_i - \si_i^* E_j = 0. 
\]
This means that the first column is the sum of the second 
and third columns in the matrix (\ref{eqn:si4}). 
Using the second formula in (\ref{eqn:xi}) repeatedly, we 
see that the matrix representation of $\si_i^{*}$ is given 
as in the first matrix of Table \ref{tab:matrix}. 
Those of $\si_j^{*}$ and $\si_k^{*}$ are obtained in the 
same manner. 
Applying Lemma \ref{lem:compo2} to these results yields 
the desired representation for $c^*$ as in the last 
matrix of Table \ref{tab:matrix}. 
Now it is easy to calculate the characteristic 
polynomial of $c^*$ as in (\ref{eqn:char}). 
The assertion for its roots, namely, for the eigenvalues 
of $c^*$ is straightforward. \hfill $\Box$ \par\medskip
We recall some more rudiments from \cite{DF}. 
Given a bimeromorphic map $f$ of a compact K\"ahler surface 
$S$, there is the concept of {\sl first dynamical degree} 
$\l_1(f)$ defined by 
\[
\l_1(f) := \lim_{n\to\infty} || (f^n)^* ||^{1/n}, 
\]
where $||\cdot||$ is an operator norm on 
$\End \, H^{1,1}(S)$. 
The limit certainly exists and one has $\l_1(f) \ge 1$. 
It is usually difficult to evaluate this number 
in a simple mean. 
But there is a distinguished class of bimeromorphic maps 
for which the first dynamical degree can be equated 
to a more tractable quantity, namely, the class of maps 
which are called analytically stable. 
Here a bimeromorphic map $f : S \carl$ is said to be 
{\sl analytically stable} (AS for short) if for any $n \in \N$ 
there is no curve $V \subset S$ such that $f^n(V) \subset I(f)$. 
From \cite{DF} we have the following lemma.
\begin{lemma} \label{lem:as}
If $f : S \carl$ is an AS bimeromorphic map, then 
the first dynamical degree $\l_1(f)$ is equal to the 
spectral radius $\rho (f^*)$ of the linear map 
$f^* : H^{1,1}(S) \carl$. 
\end{lemma}
\par 
With this lemma in hand we continue to investigate the 
Coxeter transformation (\ref{eqn:coxeter}). 
\begin{proposition} \label{prop:asc} 
The birational map $c = \si_i \ci \si_j \ci \si_k$ 
enjoys the following properties: 
\begin{enumerate} 
\item its indeterminacy set is given by $I(c) = \{p_k\}$, 
\item its exceptional set is given by $\E(c) = L$ with 
image $c(\E(c)) = \{p_i\}$, 
\item its tangent map $(dc)_{p_i}$ at $p_i$ is zero, that 
is, $p_i$ is a superattracting fixed point, 
\item it is AS, and 
\item its first dynamical degree is given by $\l_1(c) = 
2 + \sqrt{5}$. 
\end{enumerate} 
\end{proposition}
{\it Proof}.  Lemma \ref{lem:L} implies that 
$\si_k^{-1}(I(\si_j)) = \{p_k\}$ and 
$\si_k^{-1} \ci \si_j^{-1}(I(\si_i)) = 
\si_k^{-1}(\{p_j\}) = \{p_k\}$. 
Thus the indeterminacy set of $c$ is given by 
$I(c) = \{p_k\}$, which proves assertion (1). 
In order to see assertion (2), we again apply 
Lemma \ref{lem:L} to obtain 
\[
\begin{array}{rclclclcl}
c(L_i) &=& \si_i \ci \si_j \ci \si_k (L_i) 
&=& \si_i \ci \si_j (L_i) 
&=& \si_i (L_i) &=& \{p_i \}, \\[2mm]
c(L_j) &=& \si_i \ci \si_j \ci \si_k (L_j) 
&=& \si_i \ci \si_j (L_j) 
&=& \si_i(\{p_j\}) &=& \{p_i\}, \\[2mm]
c(L_k) &=& \si_i \ci \si_j \ci \si_k (L_k) 
&=& \si_i \ci \si_j (\{p_k\}) 
&=& \si_i(\{p_j\}) &=& \{p_i\}. 
\end{array}
\]
This means that $\E(c)$ is given by the union 
$L = L_i \cup L_j \cup L_k$ with $c(\E(c)) = \{p_i\}$. 
Thus assertion (2) follows. 
From assertion (2) we notice that $p_i$ is a fixed point 
of $c$ and all points on $L_j - \{p_k\}$ and on 
$L_k-\{p_j\}$ are taken to the point $p_i$ by $c$. 
So the tangent map $(dc)_{p_i}$ is zero along the 
linearly independent directions of the lines $L_j$ 
and $L_k$ with origin at $p_i$ and hence $(dc)_{p_i}$ 
itself is zero, which proves assertion (3).  
We show assertion (4) by contradiction. 
Assume that $V \subset \ol{\Sol}(\th)$ is an 
irreducible curve such that $c^n(V) \subset 
I(c) = \{p_k\}$ for some $n \in \N$. 
If $V$ intersects the affine surface $\Sol(\th)$, 
then it cannot happen that $c^n(V) \subset \{p_k\}$, 
because $c$ is bijective on $\Sol(\th)$. 
Thus $V$ must lie in $L$ and hence 
$V = L_i$, $L_j$, or $L_k$. 
But also in this case assertion (2) implies 
that $c^n(L_a) = \{p_k\}$ for $a = i,j,k$, leading 
to a contradiction. 
Hence assertion (4) is proved. 
Finally, since $c$ is AS, 
Proposition \ref{prop:in2} and Lemma \ref{lem:as} 
immediately imply assertion (5). 
\hfill $\Box$ \par\medskip
We conjecture that the topological entropy of $c$ agrees with 
the logarithm of its first dynamical degree: 
\[
h_{\mathrm{top}}(c) = \log \l_1(c) = \log(2 + \sqrt{5}). 
\]
\section{Lefschetz Fixed Point Formula} \label{sec:lefschetz} 
We are interested in the periodic points of the 
Coxeter transformation $c : \ol{\Sol}(\th) \carl$. 
Given any $N \in \N$ we can consider the set of periodic 
points of period $N$ on the projective surface $\ol{\Sol}(\th)$, 
\[
\ol{\Per}_N(c) := \{\, X \in \ol{\Sol}(\th) - I(c^N) \,:\, 
c^N(X) = X \, \}. 
\]
as well as the set of periodic points of period $N$ on the 
affine surface $\Sol(\th)$, 
\[
\Per_N(c) := \{\, x \in \Sol(\th) \,:\, c^N(x) = x \, \}, 
\]
\par
On the other hand we have defined in (\ref{eqn:per}) the set 
$\Per_N(\k)$ of periodic points of period $N$ for the 
Pochhammer-Poincar\'e map $\wp_*$. 
By Lemma \ref{lem:coxeter}, $\Per_N(\k)$ is bijectively mapped 
onto $\Per_N(c^2) = \Per_{2N}(c)$ by the Riemann-Hilbert 
correspondence (\ref{eqn:RHtk2}) and hence 
\begin{equation} \label{eqn:perper}
\mathrm{\#}\,\Per_N(\k) = \mathrm{\#}\,\Per_{2N}(c). 
\end{equation}
Thus the main aim of this article, that is, 
the enumeration of the set $\Per_N(\k)$ is reduced to that 
of $\Per_N(c)$.  
So what we should do from now on is the following: 
\begin{itemize}
\item to count the cardinality of $\ol{\Per}_N(c)$,  
\item to relate the cardinality of $\ol{\Per}_N(c)$ 
with that of $\Per_N(c)$. 
\end{itemize}
The first task will be done with the help of Lefschetz fixed 
point formula and the second task will be by a careful 
inspection of the behavior of $c$ around $L$. 
In order to apply the Lefschetz fixed point formula, 
we first need to verify the following lemma. 
\begin{lemma} \label{lem:percurve} 
For any $N \in \N$, the Coxeter transformation 
$c : \ol{\Sol}(\th) \carl$ admits no curves of periodic points 
of period $N$. 
\end{lemma}
{\it Proof}. 
By Proposition \ref{prop:in2} the Coxeter transformation 
$c^* : H^2(\ol{\Sol}(\th),\Z) \carl$ has eigenvalues 
$0$, $-1$, $2 \pm \sqrt{5}$, among which $0$ and $2 \pm \sqrt{5}$ 
are simple eigenvalues, while $-1$ is a quadruple eigenvalue 
whose eigenspace is spanned by four eigenvectors 
\[
V_0 = 2 E_0-E_i-E_j-E_k-E_l-E_m-E_n, \quad 
V_i = E_i-E_j, \quad V_j = E_k-E_l, \quad 
V_k = E_m-E_n. 
\]
In view of (\ref{eqn:lineinf}), (\ref{eqn:in1}) and 
(\ref{eqn:FG}), there are orthogonality relations 
\begin{equation} \label{eqn:VL}
(V_a, L_b) = 0 \qquad (a = 0, i, j, k, \,\,  b = i, j, k).  
\end{equation} 
We prove the lemma by contradiction. 
Assume that $c$ admits a curve (an effective divisor) 
$D \subset \ol{\Sol}(\th)$ of periodic points of some 
period $N$. 
Then we have $(c^*)^N D = (c^N)^* D = D$ in 
$H^2(\ol{\Sol}(\th),\Z)$, where $(c^*)^N = (c^N)^*$ follows 
from the fact that $c$ is AS. 
So $(c^*)^N$ has eigenvalue $1$ with eigenvector $D$. 
This eigenvalue arises as the $N$-th power of eigenvalue 
$-1$ of $c^*$ so that $N$ must be even and $D$ must be 
a linear combination of $V_0$, $V_i$, $V_j$, $V_k$. 
Hence (\ref{eqn:VL}) implies that 
\begin{equation} \label{eqn:CLa}
(D, L_a) = 0 \qquad (a = i, j, k). 
\end{equation}
We now write $D = D' + m_i L_i + m_j L_j + m_k L_k$, 
where $D'$ is either empty or an effective divisor 
not containing $L_i$, $L_j$, $L_k$ as an irreducible 
component of it and $m_i$, $m_j$, $m_k$ are nonnegative 
integers. 
Since $(L_a, L_b) = -1$ for $a = b$ and 
$(L_a, L_b) = 1$ for $a \neq b$, 
the formula (\ref{eqn:CLa}) yields
\[
\begin{array}{rclcl}
0 &=& (D, L_i) &=& (D', L_i)-m_i+m_j+m_k, \\[1mm]
0 &=& (D, L_j) &=& (D', L_j)+m_i-m_j+m_k, \\[1mm]
0 &=& (D, L_k) &=& (D', L_k)+m_i+m_j-m_k, 
\end{array}
\]
which sum up to 
\begin{equation} \label{eqn:sum}
(D', L_i)+(D', L_j)+(D', L_k)+m_i+m_j+m_k = 0.
\end{equation} 
Since none of $L_i$, $L_j$, $L_k$ is an irreducible 
component of $D'$, the intersection number 
$(D', L_a)$ must be nonnegative for any $a = i, j, k$. 
Since $m_i$, $m_j$, $m_k$ are also nonnegative, 
formula (\ref{eqn:sum}) implies that 
$(D',L_i) = (D',L_j)= (D',L_k) =0$ 
and $m_i = m_j = m_k = 0$. 
Hence $D = D'$ and $(D,L_i) = (D,L_j) = (D,L_k) =0$. 
It follows that $D$ is an effective divisor with 
$(D, L_a) = 0$ not containing $L_a$ as its 
irreducible component for every $a = i,j,k$. 
This means that the compact curve $D$ does not intersect 
$L = L_i \cup L_j \cup L_k$ and hence must lie in 
the affine cubic surface $\Sol(\th) = \ol{\Sol}(\th)-L$. 
But no compact curve can lie in any affine variety. 
This contradiction establishes the lemma. 
\hfill $\Box$ \par\medskip 
For each $N \in \Z$ let $\vG_N \subset 
\ol{\Sol}(\th) \times \ol{\Sol}(\th)$ be the graph of the 
$N$-th iterate $c^N : \ol{\Sol}(\th) \carl$, and 
$\vD \subset \ol{\Sol}(\th) \times \ol{\Sol}(\th)$ be 
the diagonal. 
Note that $\vG_N = \vG_{-N}^{\vee}$, where 
$\vG_{-N}^{\vee}$ is the reflection of 
$\vG_{-N}$ with respect to the diagonal $\vD$. 
Moreover let $I_N \subset \ol{\Sol}(\th)$ denote the 
indeterminacy set of $c^N$. 
Then the Lefschetz fixed point formula consists of 
two equations concerning the intersection number 
$(\vG_N, \vD)$ of $\vG_N$ and $\vD$ in 
$\ol{\Sol}(\th) \times \ol{\Sol}(\th)$,  
\begin{eqnarray} 
(\vG_N, \vD) &=& \sum_{q=0}^4 (-1)^q \, \Tr\,[\, 
(c^N)^* : H^q(\ol{\Sol}(\th),\Z) \carl\,], 
\label{eqn:lfpf1} \\
(\vG_N, \vD) &=& \mathrm{\#} \, \ol{\Per}_N(c) 
+ \sum_{p \in I_N} \mu((p,p), \vG_N \cap \vD), 
\label{eqn:lfpf2} 
\end{eqnarray} 
where $\mu((p,p), \vG_N \cap \vD)$ denotes the multiplicity 
of intersection between $\vG_N$ and $\vD$ at $(p,p)$. 
Lemma \ref{lem:percurve} guarantees that all terms involved in 
(\ref{eqn:lfpf1}) and (\ref{eqn:lfpf2}) are well defined and 
finite. 
\begin{lemma} \label{lem:lfpf1} 
Formula $(\ref{eqn:lfpf1})$ becomes 
$(\vG_N, \vD) = (2+\sqrt{5})^N + (2-\sqrt{5})^N + 4 (-1)^N + 2$. 
\end{lemma} 
{\it Proof}. 
We put 
$T_N^q = \Tr\,[\, (c^N)^* : H^q(\ol{\Sol}(\th),\Z) \carl\,]$. 
Because $\ol{\Sol}(\th)$ is a smooth rational surface, 
\[
H^q(\ol{\Sol}(\th), \Z) \cong \left\{ 
\begin{array}{cl}
\Z \qquad & (q = 0,4), \\[1mm] 
0  \qquad & (q = 1,3). 
\end{array}
\right. 
\]
Naturally we have $T_N^0 = 1$ and $T_N^1 = T_N^3 = 0$. 
Since $c$ and so $c^N$ are birational, we have $T_N^4 = 1$. 
By assertion (4) of Proposition \ref{prop:asc} the map 
$c$ is AS, and so Lemma \ref{lem:compo} implies that 
$(c^N)^* = (c^*)^N : H^2(\ol{\Sol}(\th),\Z) \carl$. 
Recall that $c^*$ has eigenvalues $0$, $-1$ and 
$2 \pm \sqrt{5}$, where the eigenvalue $-1$ is 
quadruple while the remaining ones are simple 
(see Proposition \ref{prop:in2}). 
Thus we have 
$T_N^2 = 0^N + 4(-1)^N + (2+\sqrt{5})^N + (2-\sqrt{5})^N$. 
Substituting these data into (\ref{eqn:lfpf1}) yields 
the assertion of the lemma. \hfill $\Box$ \par
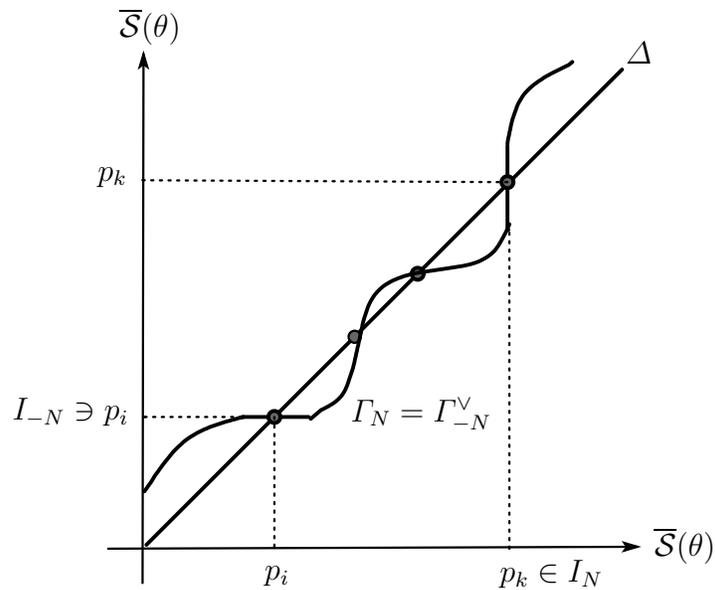
\begin{figure}[t]
\begin{center}
\unitlength 0.1in
\begin{picture}(33.60,29.90)(1.90,-30.80)
%
\special{pn 13}%
\special{pa 690 2900}%
\special{pa 3460 2890}%
\special{fp}%
\special{sh 1}%
\special{pa 3460 2890}%
\special{pa 3393 2870}%
\special{pa 3407 2890}%
\special{pa 3393 2910}%
\special{pa 3460 2890}%
\special{fp}%
%
\special{pn 13}%
\special{pa 870 3080}%
\special{pa 876 310}%
\special{fp}%
\special{sh 1}%
\special{pa 876 310}%
\special{pa 856 377}%
\special{pa 876 363}%
\special{pa 896 377}%
\special{pa 876 310}%
\special{fp}%
%
\special{pn 20}%
\special{pa 890 2880}%
\special{pa 3380 390}%
\special{fp}%
\put(35.5000,-29.8000){\makebox(0,0)[lb]{$\ol{\Sol}(\th)$}}%
\put(7.5000,-2.6000){\makebox(0,0)[lb]{$\ol{\Sol}(\th)$}}%
\put(34.0000,-3.6000){\makebox(0,0)[lb]{$\vD$}}%
\put(27.4000,-29.6000){\makebox(0,0)[lt]{$p_k\in I_N$}}%
%
\special{pn 20}%
\special{pa 2780 790}%
\special{pa 2780 1230}%
\special{fp}%
\special{pa 2780 1230}%
\special{pa 2780 1230}%
\special{fp}%
%
\special{pn 20}%
\special{sh 0.600}%
\special{ar 2780 980 32 32  0.0000000 6.2831853}%
%
\special{pn 20}%
\special{pa 1390 2210}%
\special{pa 1740 2210}%
\special{fp}%
%
\special{pn 20}%
\special{sh 0.600}%
\special{ar 1560 2210 32 32  0.0000000 6.2831853}%
%
\special{pn 13}%
\special{pa 2790 1200}%
\special{pa 2790 2900}%
\special{dt 0.045}%
\special{pa 2790 2900}%
\special{pa 2790 2899}%
\special{dt 0.045}%
%
\special{pn 13}%
\special{pa 1390 2210}%
\special{pa 880 2210}%
\special{dt 0.045}%
\special{pa 880 2210}%
\special{pa 881 2210}%
\special{dt 0.045}%
%
\special{pn 20}%
\special{sh 0.600}%
\special{ar 2310 1460 32 32  0.0000000 6.2831853}%
%
\special{pn 13}%
\special{sh 0.600}%
\special{ar 1980 1790 32 32  0.0000000 6.2831853}%
%
\special{pn 20}%
\special{pa 2790 1200}%
\special{pa 2776 1230}%
\special{pa 2761 1260}%
\special{pa 2745 1289}%
\special{pa 2728 1315}%
\special{pa 2710 1340}%
\special{pa 2688 1361}%
\special{pa 2665 1379}%
\special{pa 2640 1394}%
\special{pa 2612 1407}%
\special{pa 2583 1417}%
\special{pa 2553 1426}%
\special{pa 2521 1434}%
\special{pa 2489 1440}%
\special{pa 2455 1445}%
\special{pa 2421 1450}%
\special{pa 2387 1455}%
\special{pa 2352 1461}%
\special{pa 2318 1467}%
\special{pa 2284 1473}%
\special{pa 2250 1482}%
\special{pa 2218 1491}%
\special{pa 2186 1503}%
\special{pa 2156 1516}%
\special{pa 2129 1531}%
\special{pa 2103 1549}%
\special{pa 2080 1569}%
\special{pa 2060 1592}%
\special{pa 2043 1618}%
\special{pa 2028 1646}%
\special{pa 2015 1676}%
\special{pa 2005 1706}%
\special{pa 1995 1738}%
\special{pa 1987 1771}%
\special{pa 1979 1803}%
\special{pa 1972 1835}%
\special{pa 1965 1866}%
\special{pa 1958 1897}%
\special{pa 1950 1928}%
\special{pa 1942 1959}%
\special{pa 1934 1989}%
\special{pa 1924 2020}%
\special{pa 1914 2050}%
\special{pa 1902 2080}%
\special{pa 1888 2110}%
\special{pa 1871 2137}%
\special{pa 1850 2160}%
\special{pa 1825 2180}%
\special{pa 1798 2196}%
\special{pa 1769 2211}%
\special{pa 1750 2220}%
\special{sp}%
%
\special{pn 20}%
\special{pa 1400 2210}%
\special{pa 1368 2219}%
\special{pa 1337 2228}%
\special{pa 1305 2237}%
\special{pa 1274 2248}%
\special{pa 1244 2258}%
\special{pa 1214 2270}%
\special{pa 1185 2283}%
\special{pa 1157 2297}%
\special{pa 1130 2313}%
\special{pa 1105 2331}%
\special{pa 1080 2350}%
\special{pa 1057 2371}%
\special{pa 1035 2393}%
\special{pa 1013 2416}%
\special{pa 993 2441}%
\special{pa 973 2466}%
\special{pa 953 2493}%
\special{pa 934 2519}%
\special{pa 915 2547}%
\special{pa 897 2574}%
\special{pa 880 2600}%
\special{sp}%
%
\special{pn 20}%
\special{pa 2780 780}%
\special{pa 2787 747}%
\special{pa 2794 714}%
\special{pa 2801 682}%
\special{pa 2810 650}%
\special{pa 2820 619}%
\special{pa 2831 590}%
\special{pa 2844 561}%
\special{pa 2860 535}%
\special{pa 2877 510}%
\special{pa 2898 488}%
\special{pa 2920 467}%
\special{pa 2945 447}%
\special{pa 2971 429}%
\special{pa 2999 412}%
\special{pa 3028 396}%
\special{pa 3058 381}%
\special{pa 3088 366}%
\special{pa 3119 351}%
\special{pa 3120 350}%
\special{sp}%
%
\special{pn 13}%
\special{pa 2780 980}%
\special{pa 870 970}%
\special{dt 0.045}%
\special{pa 870 970}%
\special{pa 871 970}%
\special{dt 0.045}%
%
\special{pn 13}%
\special{pa 1560 2210}%
\special{pa 1560 2890}%
\special{dt 0.045}%
\special{pa 1560 2890}%
\special{pa 1560 2889}%
\special{dt 0.045}%
\put(7.9000,-9.0000){\makebox(0,0)[rt]{$p_k$}}%
\put(1.9000,-21.1000){\makebox(0,0)[lt]{$I_{-N}\ni p_i$}}%
\put(15.1000,-29.9000){\makebox(0,0)[lt]{$p_i$}}%
\put(19.6000,-21.1000){\makebox(0,0)[lt]{$\vG_N = \vG_{-N}^{\vee}$}}%
\end{picture}%
\end{center}
\caption{The indeterminacy point $p_k$ of $c^N$ 
is a superattracting fixed point of $c^{-N}$}
\label{fig:cubic5}
\end{figure}

\medskip 
\begin{lemma} \label{lem:lfpf2} 
Formula $(\ref{eqn:lfpf2})$ becomes 
$(\vG_N, \vD) = \mathrm{\#} \, \ol{\Per}_N(c) + 1$ with 
$\mathrm{\#} \, \ol{\Per}_N(c) = 
\mathrm{\#} \, \Per_N(c) + 1$.
\end{lemma} 
{\it Proof}. 
By Proposition \ref{prop:asc}, for any $N \in \N$, 
the point $p_k$ is the unique indeterminacy point of 
$c^N$ and the point $p_i$ is the unique fixed point of 
$c^N$ on $L$. 
Namely we have 
$I_N = \{p_k\}$ and $\ol{\Per}_N(c) = \Per_N(c) \cup \{p_i\}$,  
which implies that formula (\ref{eqn:lfpf2}) is rewritten as 
\begin{equation} \label{eqn:lfpf2a}
\begin{array}{rcl} 
(\vG_N, \vD) &=& \mathrm{\#} \, \ol{\Per}_N(c) 
+ \mu((p_k,p_k), \vG_N \cap \vD),   \\[2mm]  
\mathrm{\#} \, \ol{\Per}_N(c) &=& 
\mathrm{\#} \, \Per_N(c) + \nu(p_i, c^N), 
\end{array}
\end{equation} 
where $\nu(p_i, c^N)$ is the local index of the map 
$c^N$ around the fixed point $p_i$. 
By assertion (3) of Proposition \ref{prop:asc}, for any 
$N \in \N$, the point $p_i$ is a superattracting fixed point 
of $c^N$ and so $\det(I- (dc^N)_{p_i}) = \det(I-O) = 1$. 
This means that $\nu(p_i, c^N) = 1$. 
Likewise, since $p_k$ is a superattracting fixed 
point of $c^{-N} = (c^{-1})^N$ where $c^{-1} = 
\si_k \circ \si_j \circ \si_i$ (see Figure \ref{fig:cubic5}), 
the same reasoning as above with $c$ replaced by $c^{-1}$ 
yields $\nu(p_k, c^{-N}) = 1$. 
Therefore we have  
\[
\mu((p_k,p_k), \vG_N \cap \vD) = 
\mu((p_k,p_k), \vG_{-N}^{\vee} \cap \vD ) = 
\mu((p_k,p_k), \vG_{-N} \cap \vD) = \nu(p_k, c^{-N}) = 1. 
\]
These arguments imply that (\ref{eqn:lfpf2a}) is equivalent 
to the statement of the lemma. 
\hfill $\Box$ \par\medskip
Putting Lemmas \ref{lem:lfpf1} and \ref{lem:lfpf2} together, 
we have established the following theorem. 
\begin{theorem} \label{thm:periodic} 
For any $N \in \N$, the cardinalities of periodic points 
of period $N$ are give by 
\begin{equation} \label{eqn:carperiod}
\begin{array}{rcl}
\mathrm{\#} \, \ol{\Per}_N(c) &=& 
(2+\sqrt{5})^N + (2-\sqrt{5})^N + 4 (-1)^N + 1, \\[2mm]
\mathrm{\#} \, \Per_N(c) &=& 
(2+\sqrt{5})^N + (2-\sqrt{5})^N + 4 (-1)^N . 
\end{array} 
\end{equation}
\end{theorem}
Then our main theorem (Theorem \ref{thm:main}) is an 
immediate consequence of (\ref{eqn:perper}) and the second 
formula of (\ref{eqn:carperiod}). 
Thus the proof of Theorem \ref{thm:main} has just been 
completed. 
\par 
In this article we have seen that geometry of cubic surfaces 
and dynamics on them played an important part in understanding 
an aspect of the global structure of the sixth Painlev\'e equation. 
Their relevance to other aspects will be explored elsewhere. 

\end{document}